\documentclass[preprint,12pt]{elsarticle}

\usepackage{amssymb}
\usepackage{amsmath}
\usepackage{siunitx} %per unità di misura
\usepackage{lineno}

\usepackage{natbib}
\biboptions{sort&compress}
\usepackage{booktabs}
\usepackage{xcolor}
\usepackage{ulem}
\usepackage{caption}
\usepackage{subcaption}
\journal{Computational Physics}

\begin{document}

\begin{frontmatter}



\title{Well-balanced high-order method for non-conservative hyperbolic PDEs with source terms: application to one-dimensional blood flow equations with gravity}

\author[DICAM]{Chiara Colombo\corref{cor1}\fnref{MAT}} %% Author name
\ead{chiara.colombo-1@unitn.it}
\cortext[cor1]{Corresponding author}
\author[MAT]{Caterina Dalmaso}
\ead{caterina.dalmaso@unitn.it}
\author[MAT]{Lucas O. Müller}
\ead{lucas.muller@unitn.it}
\author[DICAM]{Annunziato Siviglia}
\ead{annunziato.siviglia@unitn.it}

\affiliation[DICAM]{organization={Department of Civil, Environmental and Mechanical Engineering, University of Trento},
            addressline={Via Mesiano 77}, 
            city={Trento},
            postcode={38123}, 
            country={Italy}}

\affiliation[MAT]{organization={Laboratory of Mathematics for Biology and Medicine, Department of Mathematics, University of Trento},
            addressline={Via Sommarive 14},
            city={Trento},
            postcode={38123},
            country={Italy}}

%% Abstract
\begin{abstract}
% \TODO{Abstract text: You are required to provide a concise and factual abstract which does not exceed 250 words. The abstract should briefly state the purpose of your research, principal results and major conclusions.}
The present work proposes a well-balanced finite volume-type numerical method for the solution of non-conservative hyperbolic partial differential equations (PDEs) with source terms.
The method is characterized, first, by the use of a recently introduced high-order spatial reconstruction, based on generalized Riemann problem information from the previous time level. Such reconstruction is well-balanced up to order three, compact, efficient and easy to implement.
Second, the method incorporates a well-balanced space-time evolution operator, which allows for well-balanced fully explicit time evolution.
The accuracy and efficiency of the method are assessed on both a scalar problem (Burgers' equation) and a nonlinear PDE system (hyperbolized one-dimensional blood flow equations with gravity and friction, and with variable mechanical and geometrical properties). The well-balanced property is verified by showing that numerically-determined stationary solutions are preserved up to machine precision. The order of accuracy in space and time is validated through empirical convergence rate studies. Additionally, the performance of the method is assessed on a network of 86 arteries, under both stationary and transient conditions.
\end{abstract}

\begin{graphicalabstract}

\vspace{1cm}
\includegraphics[width=0.95\linewidth]{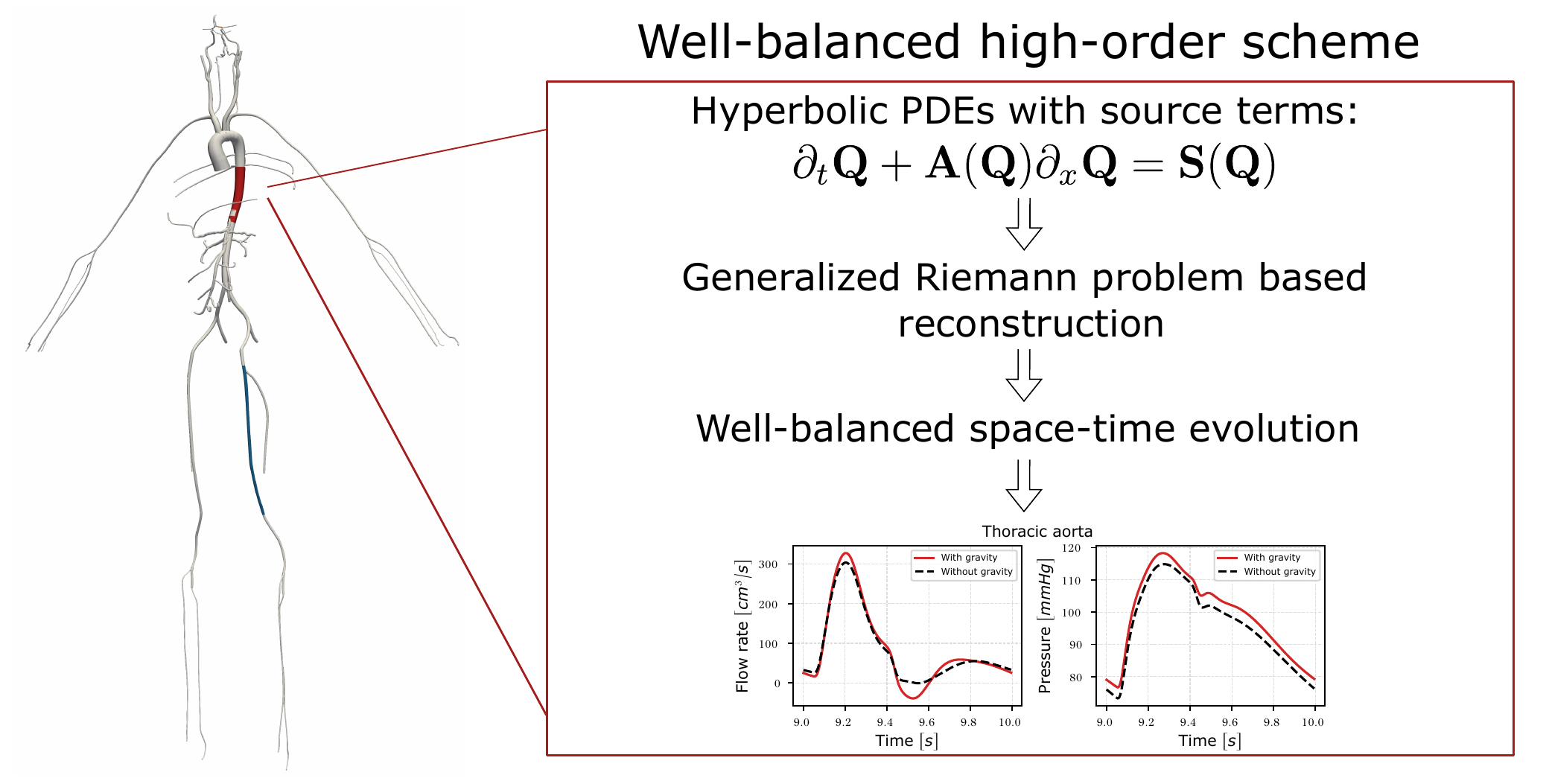}

\end{graphicalabstract}

\begin{highlights}
%\item Highlights should consist of 3 to 5 bullet points, each a maximum of 85 characters, including spaces (guardando gli altri lavori non mi sembrano così fiscali a riguardo).
\item We propose a well-balanced high-order method for non-conservative hyperbolic PDEs.
\item The method uses a GRP-based reconstruction procedure, well-balanced up to order three.
\item The method incorporates a well-balanced space-time evolution operator.
\item Steady-state solutions of Burgers' eq. and BFEs are preserved up to machine precision.
\item The method is accurate and efficient for BFEs with gravity on vascular networks.

\end{highlights}

%% Keywords
\begin{keyword}
well-balanced numerical method \sep 
high order \sep
path-conservative approach \sep
non-conservative hyperbolic PDEs \sep
GRP-based reconstruction \sep
gravity \sep
blood flow


% \TODO{TO DO: You are required to provide 1 to 7 keywords for indexing purposes. Keywords should be written in English. Please try to avoid keywords consisting of multiple words (using "and" or "of").}

\end{keyword}

\end{frontmatter}

%\linenumbers

%%%%%%%%%%%%%%%%%%%%%%%%%%%%% main text %%%%%%%%%%%%%%%%%%%%%%%%%%%%%%%%%%%%%
\section{Introduction}

Several physical phenomena can be modeled through non-conservative hyperbolic PDEs with source terms. 
Notable examples include the shallow water equations (SWEs) with moving bottom topography~\citep{introList06}, the multilayer SWEs with variable density~\citep{fernandez2022arbitrary}, the Euler equations with gravity~\citep{introList13,bermudez2016numerical},
%the description of multi-phase flows~\cite{introBase02},  %referenza non va bene. Nelle eq. non c'è termine sorgente.
and the characterization of blood flow within vessels having variable physical properties~\cite{PATHS, introBase01}.
The solution of such equations in the presence of discontinuities (shock waves) poses several challenges, as the conservative form of the equations no longer holds. One of the main issues is the definition of jump relations across a discontinuity, which has been addressed by~\citet{introNum1} interpreting the jump relations as Borel measures dependent on the path connecting the two sides of the discontinuity.

A second issue is the description of steady-state solutions in the presence of geometric- and/or algebraic-type source terms, which can be discontinuous~\cite{introNum5}. 
Such solutions should be accurately captured by the numerical methods used to solve the considered PDEs. Moreover, they should be recovered after small perturbations in the initial data, in order to avoid spurious oscillations that might cause significant deviations from the exact solution~\cite{introNum2, GF}.
A numerical scheme able to preserve the correct balance between advective and source terms was initially introduced in 1994 for the SWEs by~\citet{introNum3}, paving the way for the numerous first- and high-order approaches developed in the last thirty years for different systems of balance laws~\cite{introNum2, GF, introList01, introList02, introList03, introList04, introList05, introNum4}, with a particular focus on variants of the SWEs~\cite{introList07, introList08, introList09, introList10, fernandez2022arbitrary} and Euler equations~\cite{introList11, introList12, introList13, introList14, introList15,bermudez2016numerical}.
% \TODO{ref 1-34 in \cite{GF}, and  R4, R6, R7}.
Nowadays, numerous methods are also available for the solution of blood flow equations (BFEs), tailored to the preservation of some or all stationary solutions of the system (well-balanced and fully well-balanced methods, respectively). Many of these methods were devised to treat space-varying geometrical and physical properties of the vasculature, friction and gravity~\cite{PATHS, delestre2013well, introlist17, li2018well, introlist16, britton2020well, spilimbergo2021one,  MurilloG}. 
% \TODO{Include overview of references in https://www.sciencedirect.com/science/article/pii/S0021999122009329}  

Recently, \citet{introNum4} and~\citet{GF} proposed a generic strategy for the development of well-balanced high-order methods in the finite volume and discontinuous Galerkin frameworks, testing them on hyperbolic PDEs with continuous and discontinuous source terms. 
One of the main obstacles to the development of well-balanced methods is the construction of a well-balanced reconstruction operator~\cite{introNum4}, since standard reconstruction operators such as ENO, WENO and MUSCL~\cite{ENO, WENO, MUSCL} are not usually well-balanced. Indeed, these methods are based on standard interpolation techniques, and there is no guarantee that stationary solutions of the considered PDEs belong to the same class of reconstruction functions. A strategy to modify standard reconstruction operators so that they are well-balanced for all stationary solutions of the considered PDEs has been proposed by~\citet{CastroRec}. However, this approach
% , whilst widely used in the literature \TODO{citation?}, 
requires the solution of a nonlinear system of equations, which might have multiple or no solutions, and might be computationally expensive~\cite{introNum4}. 

% As an example, the standard formulation of the WENO reconstruction~\cite{WENO} is not well-balanced for stationary solutions of the BFEs with variable geometrical properties~\cite{PATHS}, and needs to be modified e.g. as proposed by~\citet{muller2013well}.

In this paper, we present a numerical method that combines the novel centered reconstruction technique by~\citet{GRP}
with the well-balancing approach by~\citet{introNum4} and~\citet{GF}. 
The method is developed within the path-conservative framework, using the high-order solver of~\citet{LTS} and the well-balancing strategy of~\citet{GF}, which relies on the knowledge of families of stationary solutions.
The reconstruction technique~\citep{GRP} requires, at each time step and in every computational cell, the cell average of the solution at the current time level and the solution of a generalized Riemann problem (GRP) at the cell interfaces from the previous time level. Provided that the GRP solver being used is well balanced, this minimal data dependence guarantees the well-balancing of the reconstruction technique by construction up to the third order, while also making the technique highly efficient and easy to implement.

In the following, we show that the resulting numerical method is well-balanced up to order three, testing it on both a scalar problem represented by the Burgers' equation, and a PDE system given by the hyperbolized one-dimensional BFEs presented in~\cite{hypBFEs}, including variable geometrical and mechanical vessel properties, friction and a spatially variable $C^0$ gravity term. The gravity term, in particular, is defined, depending on the considered test, either as a constant quantity, a smooth function or a polyline computed from anatomical data.
%obtained as the projection of the gravitational acceleration on the vessel axis.

% Currently available 1D blood flow models with gravity are characterized by significant simplifying assumptions. 
% \citet{MurilloG} propose a model able to describe a head-up tilt under the assumption of coplanar vessels, whose elevation with respect to the heart varies linearly along their axis. 
% \citet{Fois} propose a more complex set-up, where gravity along each vascular segment is defined as the projection of the gravity acceleration vector on the vessel axis. Still, all vessels are assumed to be coplanar, so gravity variability is not fully captured. 
% Consequently, our approach aims to extend these simplified characterizations, removing the assumption of vascular coplanarity, in order to obtain a greater physiological accuracy. This choice poses a marked numerical challenge. Indeed, to the best of our knowledge, no finite volume type method currently available in the literature can preserve the hydrostatic pressure solution of the blood flow equations in this complex setup.

The rest of the work is structured as follows. Section~\ref{Numerical method} illustrates the proposed well-balanced numerical method, detailing the spatial reconstruction procedure~\citep{GRP} and our modifications to the local space-time prediction and time evolution steps described in~\cite{LTS}. Section~\ref{Numerical tests} presents the numerical tests performed to assess the accuracy and efficiency of the solver. Numerical results are discussed in section~\ref{results}. Conclusions are drawn in section~\ref{Concl}.

\section{Numerical method} \label{Numerical method}

We consider the following system
% In this section, we propose a high-order well-balanced numerical scheme suitable for solving non-conservative hyperbolic systems of PDEs. The system of interest reads
\begin{linenomath}
\begin{equation}
\begin{cases}
\partial_t \mathbf{Q} + \mathbf{A}(\mathbf{Q})\partial_x \mathbf{Q} = \mathbf{S}(\mathbf{Q}), \quad x\in \Omega, \quad t\in \mathcal{T}, \\
\mathbf{Q}(x,0) = \mathbf{Q}_{0}(x),
\end{cases}
\label{eq1}
\end{equation}
\end{linenomath}
where $t\in \mathcal{T}$ and $x\in \Omega$ are time and space independent variables, with $\Omega = [x_A,x_B]\subset\mathbb{R}$ a one-dimensional spatial domain and $\mathcal{T}=[0,t^K]\subset\mathbb{R}_0^+$ a temporal domain. $\mathbf{Q}(x,t)\in\mathcal{B}_Q$ is the state vector, with $\mathcal{B}_Q\subset\mathbb{R}^v$ the space of the admissible states, $\mathbf{Q}_0(x)\in\mathcal{B}_Q$ is the initial condition, $\mathbf{S}(\mathbf{Q}(x,t))\in\mathbb{R}^v$ is the source vector, and $\mathbf{A}(\mathbf{Q}(x,t))\in\mathbb{R}^{v\times v}$ is the system matrix with real and distinct eigenvalues.
If $\mathbf{A}(\mathbf{Q})$ is the Jacobian of the system, then the previous equation reduces to a classical balance law.

Our goal is to construct a numerical scheme that preserves stationary solutions in the framework of path-conservative numerical methods. We begin by identifying the stationary solutions of problem~\eqref{eq1} as described below. 
Let $\mathbf{Q}^*(x)\in\mathcal{B}_Q$, $\forall x\in\Omega$,  be the stationary solution of problem~\eqref{eq1}, then it holds
\begin{linenomath}
\begin{equation}
\mathbf{A}(\mathbf{Q}^*)\partial_x \mathbf{Q}^* = \mathbf{S}(\mathbf{Q}^*).
\label{eq2}
\end{equation}
\end{linenomath}
Subtracting eq.~\eqref{eq2} from eq.~\eqref{eq1} \citep{GF}, we obtain an equivalent problem given by
\begin{linenomath}
\begin{equation}
\begin{cases}
\partial_t \mathbf{Q} + \mathbf{A}(\mathbf{Q})\partial_x \mathbf{Q}-\mathbf{A}(\mathbf{Q}^*)\partial_x \mathbf{Q}^* = \mathbf{S}(\mathbf{Q})-\mathbf{S}(\mathbf{Q}^*), \quad x\in \Omega, \quad t\in \mathcal{T}, \\
\mathbf{Q}(x,0) = \mathbf{Q}_{0}(x).
\end{cases}
\label{eq3}
\end{equation}
\end{linenomath}
% for $x\in \Omega$, and $t\geq0$.
Then, we discretize the spatial domain $\Omega$ into $N$ computational cells $S_i=[x_{i-\frac12}, x_{i+\frac12}]$, with $i=1,\dots,N$, and the temporal domain $\mathcal{T}$ into $K$ computational cells $T^n=[t^{n},t^{n+1}]$ with $n=0,\dots , K-1$.
Integrating eq.~\eqref{eq3} in space and time over the control volume $V_i^n = S_i\times T^n$, we obtain the following numerical scheme

%A path-conservative numerical scheme of Roe type~\cite{Pares2006} for problem~\eqref{eq3} reads
% Since we are interested in a high-order numerical scheme in the framework of path-conservative numerical methods, we discretize the spatial domain into $N$ computational cells $S_i=[x_{i-\frac12}, x_{i+\frac12}]$, with $i=0,\dots,N-1$, and the time domain into $K$ computational cells $T_i=[t^n,t^{n+1}]$ with $n=0,\dots , K-1$. Then, we %apply a numerical scheme that reads
% integrate eq.~\eqref{eq3} in space and time over the control volume $V_i = [x_{i-\frac12},x_{i+\frac12}]\times[t^n, t^{n+1}]$, and we get
\begin{linenomath}
\begin{equation}
    \mathbf{Q}_i^{n+1} = \mathbf{Q}_i^n - \frac{1}{\Delta x}\bigg(\mathbf{B}_i -\mathbf{B}_i^*\bigg)
    %\int_{t^n}^{t^{n+1}}\int_{x_{i-\frac12}}^{x_{i+\frac12}} \mathbf{A}(\mathbf{Q})\partial_x\mathbf{Q}dxdt  \\
    %&+ \frac{1}{\Delta x}\int_{t^n}^{t^{n+1}}\int_{x_{i-\frac12}}^{x_{i+\frac12}} \mathbf{A}(\mathbf{Q}^*)\partial_x\mathbf{Q}^* dxdt \\
    - \frac{\Delta t^n}{\Delta x}\bigg(\mathbf{D}_{i+\frac12}^- + \mathbf{D}_{i-\frac12}^+\bigg) + \Delta t^n \bigg(\mathbf{S}_i -\mathbf{S}^*_i\bigg),
\label{eq4}
\end{equation}
\end{linenomath}
where 
\begin{linenomath}
\begin{equation}\label{eq:CA}
\mathbf{Q}_i^n \approx \frac{1}{\Delta x} \int_{x_{i-\frac12}}^{x_{i+\frac12}} \mathbf{Q}(x,t^n)dx,
\end{equation}
\end{linenomath}
\begin{linenomath}
\begin{equation} \label{eq:NCs}
\mathbf{B}_i \approx \int_{t^n}^{t^{n+1}}\int_{x_{i-\frac12}}^{x_{i+\frac12}} \mathbf{A}(\mathbf{Q})\partial_x\mathbf{Q}dxdt, \quad
\mathbf{B}_i^* \approx \int_{t^n}^{t^{n+1}}\int_{x_{i-\frac12}}^{x_{i+\frac12}} \mathbf{A}(\mathbf{Q}^*)\partial_x\mathbf{Q}^*dxdt,
\end{equation}
\end{linenomath}
\begin{linenomath}
\begin{equation} \label{eq:Sources}
\mathbf{S}_i \approx \frac{1}{\Delta x \Delta t^n} \int_{t^n}^{t^{n+1}}\int_{x_{i-\frac12}}^{x_{i+\frac12}} \mathbf{S}(\mathbf{Q})dxdt, \quad
\mathbf{S}_i^* \approx \frac{1}{\Delta x \Delta t^n} \int_{t^n}^{t^{n+1}}\int_{x_{i-\frac12}}^{x_{i+\frac12}} \mathbf{S}(\mathbf{Q}^*)dxdt,
\end{equation}
\end{linenomath}
and
\begin{linenomath}
\begin{equation}
\mathbf{D}_{i+\frac12}^{\pm} \approx \frac{1}{\Delta t^n} \int_{t^n}^{t^{n+1}} D_{i+\frac12}^{\pm}\bigg(\mathbf{Q}_{i+\frac12}^-(t),\mathbf{Q}_{i+\frac12}^+(t),\mathbf{\Psi}\bigg(\mathbf{Q}_{i+\frac12}^-(t),\mathbf{Q}_{i+\frac12}^+(t),s\bigg)\bigg)dt.
\label{eq7}
\end{equation}
\end{linenomath}
$\Delta x=x_{i+\frac12}-x_{i-\frac12}$ is the mesh spacing, and $\Delta t^n=t^{n+1}-t^{n}$ is the $n$-th time step. 
$\mathbf{Q}_i^n$ denotes the $i$-th cell average at time $t^n$, $\mathbf{B}_i$ and $\mathbf{B}_i^*$ are the non-conservative products of the $i$-th cell, $\mathbf{S}_i$ and $\mathbf{S}_i^*$ are the numerical source terms of the $i$-th cell, while $\mathbf{D}_{i+\frac{1}{2}}^\pm$ are the numerical fluctuations across the interface $x_{i+\frac{1}{2}}$. 
$D_{i+\frac12}^{\pm}\big(\mathbf{Q}_{i+\frac12}^-,\mathbf{Q}_{i+\frac12}^+,\mathbf{\Psi}\big)$ are the jump terms on the cells boundaries, also called fluctuations, that
arise from the solution of a Riemann problem (RP) at the cell interface $x_{i+\frac{1}{2}}$.
They
depend on the left and right states at the cell interface $\mathbf{Q}_{i+\frac12}^\pm(t)$, and on the integration path $\mathbf{\Psi}\big(\mathbf{Q}_{i+\frac12}^-,\mathbf{Q}_{i+\frac12}^+,s\big)$.
Here, we compute them as described in~\citep{PATHS}.
The integration path $\mathbf{\Psi}$ is a parametric arc in the parameter $s\in[0,1]$ that is used to connect the left and right states $\mathbf{Q}_{i+\frac12}^\pm(t)$. Multiple path choices are possible. Theoretical details on the choice of paths are given in~\cite{Pares2006}.
Here, we adopt the segment path for the Burgers' equation, and a modification of it for the BFEs~\cite{PATHS}.
The left and right states $\mathbf{Q}_{i+\frac12}^\pm(t)$ are space-time reconstructed data that are extrapolated to both sides of the cell interface $x_{i+\frac{1}{2}}$. 
These space-time reconstructed data  
are space-time predictions of the sought solution $\mathbf{Q}(x,t)$, 
computed with a local implicit discontinuous Galerkin scheme.
Here, we adopt a modified version of the Dumbser-Enoux-Toro GRP solver~\citep{DET} that uses spatial reconstruction polynomials as initial data.
These polynomials are computed through a reconstruction procedure that exploits GRP-based predictions from the previous time step to construct the reconstructions at the current step. Together, the GRP solver and the spatial reconstruction technique, described in detail below, provide the quantities required for evaluating the integrals in~\eqref{eq4}, thereby yielding a high-order numerical scheme in both space and time.
%Details on the spatial reconstruction technique are given in section~\ref{rec}, while the adopted GRP solver is described in section~\ref{DET}.

% Finally, the fluctuations $D_{i+\frac12}^{\pm}$ are computed as proposed in \cite{Dumbser2010}.

% Additionally, i
In order to have a fully explicit evolution of the data $\mathbf{Q}_i^{n+1}$ at each time step, we have to approximate integrals appearing in eq.~\eqref{eq4} with the desired order of accuracy.
The numerical fluctuations $\mathbf{D}_{i+\frac12}^{\pm}$~\eqref{eq7} are computed as follows. 
First, we solve a classical RP for system in~\eqref{eq3} with initial condition defined as
\begin{linenomath}
    \begin{equation}
    \mathbf{Q}(x,t^n) = 
        \begin{cases}
            \mathbf{Q}_{i+\frac{1}{2}}^- & \text{if } x<x_{i+\frac{1}{2}}, \\
            \mathbf{Q}_{i+\frac{1}{2}}^+ & \text{if } x>x_{i+\frac{1}{2}},
        \end{cases}
    \end{equation}
\end{linenomath}
in order to compute the jump terms $D^\pm_{i+\frac{1}{2}}$~\citep{PATHS}.
Then, we approximate the time integral in eq.~\eqref{eq7} using a second- or third-order quadrature rule with Gaussian quadrature nodes and Lagrange interpolation polynomials defined on these nodes.
In analogous manner, the numerical sources $\mathbf{S}_i$ and $\mathbf{S}_i^*$ in~\eqref{eq:Sources}, and the two non-conservative products $\mathbf{B}_i$ and $\mathbf{B}_i^*$ in~\eqref{eq:NCs} are approximated with the same quadrature rule in both space and time.

% To do so, we perform the following three steps at each time level, that we describe in the next subsections:
% \begin{itemize}
%     \item[1)] piecewise polynomial spatial reconstruction,
%     \item[2)] stationary solution identification,
%     \item[3)] local space-time prediction.
% \end{itemize}

% The local space-time predictions, along with the stationary solution, are then used to compute the space-time average integrals found in the source terms and non-conservative products in eq.~\eqref{eq4}. These integrals are approximated using a quadrature rule of appropriate accuracy. Here, we adopted a nodal basis function, with Gaussian quadrature nodes and Lagrange interpolation polynomials defined on these nodes as basis functions.
% Similarly, the integrals of the fluctuations~\eqref{eq7} are computed using the same quadrature rule, but employing at time $t\in[t^n,t^{n+1}]$ the solution to a classical Riemann problem with initial condition defined on 
% the space-time reconstructed data extrapolated to both sides of the cell interface.

\subsection{Stationary solution identification}
\label{SteadyState}
In this section, we describe the procedure that, at each time step $t^n$, and for each computational cell $S_i$, we use to identify a suitable stationary solution $\mathbf{Q}^*(x)$ of problem~\eqref{eq1} defined in each of the $P(=2,3)$ quadrature points $x_p$, $p=0,\dots , P-1$ used in the numerical scheme. 
To this end, we discretize each computational cell $S_i$ into $P-1$ intervals $[x_{p-1},x_{p}]$ of length $h$. Then, we identify a succession of $P$ values $\{\mathbf{Q}^*_{i,p}\}_{p=0}^{P-1}$ representing the approximated stationary solution in each node $x_p$.
Specifically, this succession is obtained by applying the Runge-Kutta method (RK) of order $P$ and with $P$ stages to eq.~\eqref{eq2}.
In this work, we consider the second-order RK with 2 stages, and the third-order RK with 3 stages. 
In both cases, the standard Butcher tableau with coefficients $a_{jk},b_j$, and $c_j$ is employed.
Therefore, we compute $\mathbf{Q}^*(x_p)\approx\mathbf{Q}^*_{i,p}$ for $p>0$ as
\begin{linenomath}
    \begin{equation} \label{eq:ssSys1}
        \begin{cases}
        \mathbf{Q}_{i,p+1}^* = \mathbf{Q}_{i,p}^* + h \sum_{j=0}^{P-1}b_jK_j, \\
        K_j = \tilde{f}(x_p + c_jh, \mathbf{Q}^*_{i,p}+h\sum_{k=0}^{P-1}a_{jk}K_k).
        \end{cases}
    \end{equation}
\end{linenomath}
Additionally, in order to identify $\mathbf{Q}_{i,0}^*$, we require that the average stationary solution in $S_i$ coincides with the cell average $\mathbf{Q}_i^n$~\cite{GF}. Particularly, we apply a Gaussian quadrature rule of order $P$, with nodes $x_p$ and weights $w_p$, to the average stationary solution integral and we enforce the equivalence at the discrete level. 
As a result, we compute $\mathbf{Q}^*_{i,0}$ by solving the following non-linear system of equations
\begin{linenomath}
    \begin{equation} \label{eq:ssSys2}
        w_0\mathbf{Q}_{i,0}^* + \sum_{p=1}^{P-1} w_p\mathbf{Q}_{i,p}^*(\mathbf{Q}_{i,0}^*) = \mathbf{Q}_i^n.
    \end{equation}
\end{linenomath}
Here, the solution to this system was found by applying the standard Newton method.

% compute it as in~\cite{GF}, solving the following system:
% \begin{linenomath}
%     \begin{equation}
%         \begin{cases}
%             \mathbf{A}(\mathbf{Q}^*_i)\partial_x \mathbf{Q}^*_i = \mathbf{S}(\mathbf{Q}^*_i), \quad \forall x \in S_i, \\
%             \frac{1}{\Delta x} \int_{x_{i-\frac12}}^{x_{i+\frac12}} \mathbf{Q}^*_i(x)dx = \frac{1}{\Delta x} \int_{x_{i-\frac12}}^{x_{i+\frac12}} \mathbf{Q}(x,t^{n})dx,  \\
%             \displaystyle \min_{\mathbf{Q}_i^*(x)\in\mathcal{B}_Q}\bigg\{\int_{x_{i-\frac12}}^{x_{i+\frac12}} \big(\mathbf{Q}^*_i(x) - \mathbf{Q}(x,t^{n})\big)^2 \bigg\}. \\
%         \end{cases}
%     \end{equation}
% \end{linenomath}
% We decided to approximate the integrals in the previous system using appropriate quadrature formulas to obtain a non-linear system of equations, that we later solved applying a Runge-Kutta method of adequate order of accuracy.

\subsection{The well-balanced DET solver}\label{DET}
%In this section, we present the procedure that we use to approximate the sought solution $\mathbf{Q}(x,t)$ inside $V_i^n$, to be used in the integrals of eq.~\eqref{eq4}, in the form of a space-time polynomial $\mathbf{Q}_i^{ST,n}(x,t)$.
In this section, we present the procedure for constructing a space-time polynomial $\mathbf{Q}_i^{ST,n}(x,t)$ that approximates the solution $\mathbf{Q}(x,t)$ within $V_i^n$ and is employed in the evaluation of the integrals in eq.~\eqref{eq4}.
Particularly, we apply a modified version of the Dumbser-Enaux-Toro (DET) GRP solver~\cite{DET} that is based on the calculation of the so-called deviations. 
The deviation $\mathbf{d}_i^n(x,t)$ is defined as the difference between the space-time polynomial $\mathbf{Q}_i^{ST,n}(x,t)$ and the stationary solution $\mathbf{Q}^*_i(x)$, namely
% Once that the reconstruction polynomials and the stationary solutions in each computational cell are known, we define the so-called deviations 
\begin{linenomath}
    \begin{equation}
        \mathbf{d}_i^n(x,t) = \mathbf{Q}_i^{ST,n}(x,t) - \mathbf{Q}^*_i(x), \quad \forall (x,t)\in V_i^n.
    \label{eq:dev}
    \end{equation}
\end{linenomath}
% The deviation $\mathbf{d}_i^n(x,t)$ expresses the distance from the stationary solution $\mathbf{Q}^*_i(x)$ of the sought solution inside $V_i$ in the form of a space-time polynomial $\mathbf{Q}_i^{ST,n}(x,t)$, that we are going to identify through the use of a modified version of the Dumbser-Enaux-Toro (DET) GRP solver~\cite{DET}.
In order to compute integrands in eq.~\eqref{eq4} at the desired time levels, we need to know the solution of the GRP for system in~\eqref{eq3} with initial condition defined as
\begin{linenomath}
    \begin{equation}
    \mathbf{Q}(x,t^n) = 
        \begin{cases}
            \mathbf{w}_{i-1}^n(x), & \text{ if }x<x_{i-\frac{1}{2}}, \\
            \mathbf{w}_{i}^n(x), & \text{ if }x>x_{i-\frac{1}{2}},
        \end{cases}
        \label{IC_GRP}
    \end{equation}
\end{linenomath}
where $\mathbf{w}^n_i(x)$ is the reconstruction polynomial of order $P-1$ in $S_i$ at time $t^n$.
The DET solver~\citep{DET} finds the solution to this GRP at the desired time points by first performing a space-time evolution procedure, and then solving classical RPs at those time points using the space-time predictions as initial conditions. 
% In particular, it uses a local space-time discontinuous Galerkin finite element scheme of order $P$ for the evolution, which returns locally evolved data in space and time, the space-time predictions $\mathbf{Q}^{ST,n}_i$.
In particular, locally evolved data in space and time, the space-time predictions $\mathbf{Q}^{ST,n}_i$, are computed through a local space-time discontinuous Galerkin finite element scheme of order $P$.
%, that are used to approximate the integrals in~\eqref{eq4}.
Given the structure of the adopted numerical scheme, these predictions $\mathbf{Q}_i^{ST,n}$ are space-time polynomials $\mathbf{Q}_i^{ST,n}(x,t)$ evaluated at spatial and temporal quadrature nodes, i.e. $x_{i-\frac12},x_{i+\frac12}$ and $t^{n},t^{n+1}$ for a second-order method, and $x_{i-\frac12}, x_i,x_{i+\frac12}$ and $t^{n},t^{n+\frac12},t^{n+1}$ for a third-order method.
The classical Riemann solver instead provides the solution of the GRP at the time points along interfaces where such solution is required for the computation of the numerical fluctuations $\mathbf{D}_{i\pm\frac12}^\pm$ ~\eqref{eq7} and the reconstruction polynomials, by solving classical RPs and using space-time predictions at the previous time level on the left and on the right of the same interface as initial data.
Before proceeding with the description of the DET solver, 
%we underline that our modified version of the DET solver is equivalent to the original version of the DET solver if we assume that the stationary solution $\mathbf{Q}^*_i(x)$ is always zero. Additionally, 
we underline that whenever we refer to a local quantity, we drop both the temporal superscript $n$ and the spatial subscript $i$, keeping only the subscript $h$ to recall that we are in the local framework, for example $\mathbf{Q}_h^{ST} = \mathbf{Q}_i^{ST,n}$.

The local framework is defined by the transformation of the computational volume $V_i^n$ into the reference space-time element $V_h = [0,1]\times[0,1]$ with the reference coordinate $\xi$ and $\tau$ as $x = x_{i-\frac12} + \Delta x \xi$, and $t = t^n + \Delta t^n \tau$.
Hence, problem~\eqref{eq3} becomes
\begin{linenomath}
\begin{equation}
\begin{cases}
\begin{aligned}
\partial_\tau \mathbf{Q}^{ST}_h 
+\frac{\Delta t^n}{\Delta x} \mathbf{A}(\mathbf{Q}^{ST}_h)\partial_\xi \mathbf{Q}^{ST}_h
-\frac{\Delta t^n}{\Delta x}\mathbf{A}(\mathbf{Q}^*_h)\partial_\xi \mathbf{Q}^*_h = \\ 
\Delta t^n \mathbf{S}(\mathbf{Q}^{ST}_h)-\Delta t^n \mathbf{S}(\mathbf{Q}^*_h),
\end{aligned} \\
\mathbf{Q}^{ST}_h(\xi, 0) = \mathbf{w}_h(\xi).
\end{cases}
\label{eq10}
\end{equation}
\end{linenomath}
We note that this change of coordinate does not affect the mesh spacing $\Delta x$ or the time step $\Delta t^n$, as both are constant quantities.
Since $\partial_\tau \mathbf{Q}^*_h = 0$, 
%due to the linearity of the temporal derivative operator, 
using relation \eqref{eq:dev},
we replace the term $\partial_\tau \mathbf{Q}^{ST}_h$ with $\partial_\tau \mathbf{d}_{h}$.
Additionally, we write the space-time predictor in terms of deviations as $\mathbf{Q}^{ST}_h(\xi,\tau) = \mathbf{d}_h(\xi,\tau) + \mathbf{Q}^{*}_h(\xi)$. As a result, we transform the unknown of our problem from $\mathbf{Q}^{ST}_h(\xi,\tau)$ to $\mathbf{d}_h(\xi,\tau)$.
We also note that if both the non-conservative product $\mathbf{A}(\mathbf{Q})\partial_\xi\mathbf{Q}$ and the source term $\mathbf{S}(\mathbf{Q})$ are linear operators, then eq.~\eqref{eq10} can be written only in terms of deviations $\mathbf{d}_h(\xi,\tau)$, without the stationary solution $\mathbf{Q}^*_h(\xi)$ explicitly appearing in the formulation.

Then, we multiply eq.~\eqref{eq10} by a space-time basis function $\theta(\xi,\tau)\in \mathbb{P}_{m,m}$ with $m$ the degree of the polynomial, and we integrate it over $V_h$. After applying integration by parts only in time, we obtain
\begin{linenomath}
    \begin{equation}
    \begin{split}
            [\theta, \mathbf{d}_h]^1 
             & - \langle \partial_\tau\theta, \mathbf{d}_h \rangle_{V_h}  
            - [\theta, \mathbf{d}_h]^0
        +\frac{\Delta t^n}{\Delta x} \langle \theta, \mathbf{A}(\mathbf{Q}^{ST}_h)\partial_\xi\mathbf{Q}^{ST}_h\rangle_{V_h} \\
         & -\frac{\Delta t^n}{\Delta x} \langle \theta, \mathbf{A}(\mathbf{Q}^{*}_h)\partial_\xi\mathbf{Q}^{*}_h\rangle_{V_h} =
             \Delta t^n \langle \theta, \mathbf{S}(\mathbf{Q}^{ST}_h)\rangle_{Vh}
            - \Delta t^n \langle \theta, \mathbf{S}(\mathbf{Q}^{*}_h)\rangle_{Vh},
    \end{split}
    \label{eq11}
    \end{equation}
\end{linenomath}
where we have used the following notation for the two scalar products of two functions $f(\xi,\tau)$ and $g(\xi,\tau)$
\begin{linenomath}
    \begin{equation}
        [f,g]^\tau = \int_0^1 f(\xi,\tau)g(\xi,\tau)d\xi, \quad \langle f,g \rangle_{V_h} = \int_0^1\int_0^1 f(\xi,\tau)g(\xi,\tau)d\xi d\tau.
    \end{equation}
\end{linenomath}

We now take the element $[\theta, \mathbf{d}_h]^0$ of eq.~\eqref{eq11} at reference time $\tau=0$ and we observe that it is completely defined by the initial condition. 
The initial condition~\citep{DET} is
\begin{linenomath}
    \begin{equation}
       \mathbf{d}_{h,0}(\xi) = \mathbf{d}_h(\xi,0) = \mathbf{w}_h(\xi) - \mathbf{Q}^*_h(\xi),
    \end{equation}
\end{linenomath}
where the stationary term $\mathbf{Q}^*_h(\xi)$ is known from previous computations.
%, and the initial condition related to the space-time predictor $\mathbf{Q}_{h}^{ST}(\xi,0)$ is given by the reconstruction polynomial $\mathbf{w}_h(\xi)$.

Later, using the same basis function $\theta(\xi,\tau)$, we approximate the different quantities in \eqref{eq11} expanding them as
\begin{linenomath}
    \begin{equation}
        \mathbf{d}_h(\xi, \tau) = \sum_{l=1}^{P^2} \theta_l \widehat{\mathbf{d}}_l, \quad
        \mathbf{Q}^*_h(\xi) = \sum_{l=1}^{P^2} \theta_l \widehat{\mathbf{Q}^*}_l,
    \end{equation}
\end{linenomath}
\begin{linenomath}
    \begin{equation}
    \mathbf{A}(\mathbf{Q}^{ST}_h)\partial_\xi\mathbf{Q}^{ST}_h(\xi, \tau) = 
    \sum_{l=1}^{P^2} \theta_l \widehat{\mathbf{A}\partial_\xi \mathbf{Q}}_l, \quad
    \mathbf{A}(\mathbf{Q}^{*}_h)\partial_\xi\mathbf{Q}^{*}_h(\xi) =
    \sum_{l=1}^{P^2} \theta_l \widehat{\mathbf{A}^*\partial_\xi \mathbf{Q}^*}_l,
    \end{equation}
\end{linenomath}
\begin{linenomath}
    \begin{equation}
        \mathbf{S}(\mathbf{Q}^{ST}_h)(\xi, \tau) = \sum_{l=1}^{P^2} \theta_l \widehat{\mathbf{S}}_l =
         \sum_{l=1}^{P^2} \theta_l \mathbf{S}(\widehat{\mathbf{Q}}_l) =
        \sum_{l=1}^{P^2} \theta_l \mathbf{S}(\widehat{\mathbf{d}}_l + \widehat{\mathbf{Q}^*}_l),
    \end{equation}
\end{linenomath}
\begin{linenomath}
    \begin{equation}
        \mathbf{S}(\mathbf{Q}^{*}_h)(\xi) = \sum_{l=1}^{P^2} \theta_l \widehat{\mathbf{S}^*}_l = \sum_{l=1}^{P^2} \theta_l \mathbf{S}(\widehat{\mathbf{Q}^*}_l).
    \end{equation}
\end{linenomath}
\sloppy %to avoid inline formulas to overflow
We observe that the term $\widehat{\mathbf{A}\partial_\xi \mathbf{Q}}_l$ is obtained by assuming that $\widehat{\mathbf{A}\partial_\xi \mathbf{Q}}_l = \mathbf{A}(\widehat{\mathbf{Q}}_l)\widehat{\partial_\xi \mathbf{Q}}_l$, and by computing the expansion coefficients of the spatial derivative as
$\langle \theta_k, \theta_l\rangle_{V_h} \widehat{\partial_\xi \mathbf{Q}}_l = \langle \theta_k, \partial_\xi  \theta_l\rangle_{V_h} \widehat{\mathbf{Q}}_l$ for $k=1,\dots,P^2$~\citep{LTS}.
Similarly, $\widehat{\mathbf{A}^*\partial_\xi \mathbf{Q}^*}_l = \mathbf{A}(\widehat{\mathbf{Q}^*}_l)\widehat{\partial_\xi \mathbf{Q}^*}_l$, and the spatial derivative is computed as 
$\langle \theta_k, \theta_l\rangle_{V_h} \widehat{\partial_\xi \mathbf{Q}^*}_l = \langle \theta_k, \partial_\xi  \theta_l\rangle_{V_h} \widehat{\mathbf{Q}^*}_l$ for $k=1,\dots,P^2$.
In general, $\forall l=1,\dots,P^2$, the expansion coefficients $\widehat{\mathbf{Q}^*}_l$, $\widehat{\mathbf{A}^*\partial_\xi \mathbf{Q}^*}_l$, and $\widehat{\mathbf{S}^*}_l$, , are known because of the knowledge of the stationary solution. 
The expansion coefficients $\widehat{\mathbf{Q}}_l$ instead can be written in terms of deviations as $\widehat{\mathbf{Q}}_l = \widehat{\mathbf{d}}_l + \widehat{\mathbf{Q}^*}_l$. 
Finally, the expansions coefficients related to the deviations $\widehat{\mathbf{d}}_l$, $\forall l=1,\dots,P^2$, which are the resulting unknown of system~\eqref{eq11}, are found applying a fixed-point iteration procedure to the following system
\begin{linenomath}
    \begin{equation}
    \begin{split}
             & \big\{[\theta_k, \theta_l]^1 
             - \langle \partial_\tau\theta_k, \theta_l \rangle_{V_h}\big\} \widehat{\mathbf{d}}_l^{m+1} 
              - \Delta t^n \langle \theta_k, \theta_l\rangle_{Vh}\mathbf{S}(\widehat{\mathbf{d}}_l^{m+1} + \widehat{\mathbf{Q}^*}_l) = \\
            [\theta_k, \psi_l]^0 & \widehat{\mathbf{d}_0}_{l} 
         -\frac{\Delta t^n}{\Delta x} \langle \theta_k, \theta_l\rangle_{V_h} \mathbf{A}(\widehat{\mathbf{d}}_l^m + \widehat{\mathbf{Q}^*}_l)\big\{\langle \theta_k, \theta_l \rangle_{V_h}\big\}^{-1} \langle \theta_k, \partial_\xi\theta_l \rangle_{V_h}(\widehat{\mathbf{d}}_l^m + \widehat{\mathbf{Q}^*}_l) \\
        +\frac{\Delta t^n}{\Delta x} & \langle \theta_k, \theta_l\rangle_{V_h} \mathbf{A}(\widehat{\mathbf{Q}^*}_l)\big\{\langle \theta_k, \theta_l \rangle_{V_h}\big\}^{-1} \langle \theta_k, \partial_\xi\theta_l \rangle_{V_h}\widehat{\mathbf{Q}^*}_l
     - \Delta t^n \langle \theta_k, \theta_l\rangle_{Vh}\mathbf{S}(\widehat{\mathbf{Q}^*}_l),
    \end{split}
    \label{eq18}
    \end{equation}
\end{linenomath}
for $k=1,\dots,P^2$, where we have expanded the initial condition as $\mathbf{d}_{h,0}(\xi)=\sum_{l=1}^{P^2} \psi_l \widehat{\mathbf{d}_0}_l = \sum_{l=1}^{P^2} \psi_l (\widehat{\mathbf{w}}_l - \widehat{\mathbf{Q}^*}_l)$, and where $m$ represents the current iteration step.
We underline that due to the non-linearity of the system we decided to evaluate the second term of the right-hand side at the known expansion coefficients $\widehat{\mathbf{d}}_l^m$ of the current iteration.
Once the expansion coefficients are known, we finally retrieve the space-time predictions $\mathbf{Q}^{ST,n}_i$ in $V_i^n$. %This quantity corresponds to $\mathbf{Q}$ in eq.~\eqref{eq4}.
We remark that the modified version of the DET solver here described is equivalent to the original version of the DET solver if we assume that the stationary solution $\mathbf{Q}^*_i(x)$ is always zero.

\subsection{Spatial reconstruction}
\label{rec}
In this section, we provide a brief overview of the GRP-based reconstruction procedure used to derive second- and third-order accurate spatial reconstruction polynomials $\mathbf{w}_i^n(x)$~\cite{GRP}. 
These polynomials serve as local reconstructed data that provide the initial condition for the GRP solver, which in turns computes high-order space-time predictions of the solution $\mathbf{Q}(x,t)$. As a result, they enable the construction of a high-order scheme in both space and time.
%Further details on this methodology, including the derivation of higher order polynomials, can be found in~\cite{GRP}.

They 
%spatial reconstruction polynomials $\mathbf{w}_i^n(x)$ 
are piecewise functions of degree $P-1$, with $P$ being the order of the numerical method, defined on each computational cell $S_i$ at time $t^n$. They are expressed in terms of the cell average of the solution at the current time step $\mathbf{Q}_i^n$, and in terms of left and right boundary states at the previous time step $\mathbf{Q}_{i-\frac12}^{+,n-1}$ and $\mathbf{Q}_{i+\frac12}^{-,n-1}$. 
These states arise from the solution of a classical RP at time $t^{n-1}+\Delta t^{n-1}$ at the interfaces $x_{i\pm\frac12}$ between two neighboring cells, and are also used to compute the integrand in the time integral of the fluctuations~\eqref{eq7}.
For the left interface $x_{i-\frac12}$, the classical RP reads
\begin{linenomath}
    \begin{equation} \label{eq:RP11}
        \begin{cases}
            \partial_t \mathbf{Q} + \mathbf{A}(\mathbf{Q})\partial_x \mathbf{Q} = 0, \quad x\in \mathbb{R}, t>t^{n-1}, \\
\mathbf{Q}(x,t^{n-1}+\Delta t^{n-1}) =
\begin{cases}
    \mathbf{Q}_{i-1}^{ST,n-1}(x_{i-\frac12},t^{n-1}+\Delta t^{n-1}) & \text{if }x<x_{i-\frac12}, \\
    \mathbf{Q}_{i}^{ST,n-1}(x_{i-\frac12},t^{n-1}+\Delta t^{n-1}) & \text{if }x>x_{i-\frac12},
\end{cases}
        \end{cases}
    \end{equation}
\end{linenomath}
where $\mathbf{Q}_{i-1}^{ST,n-1}$ and $ \mathbf{Q}_{i}^{ST,n-1}$ are space-time predictions at time $t^{n-1}+\Delta t^{n-1}$ and related to cells $S_{i-1}$ and $S_i$, respectively. Their computation is described in the previous section. 
The solution to RP~\eqref{eq:RP11} can be obtained using various methods~\cite{Toro}. Here we applied the exact Riemann solver for the Burgers' equation, and a two-rarefaction approximate Riemann solver for the BFEs, which yielded to a non-linear system in the left boundary states unknowns $\mathbf{Q}_{i-\frac12}^{\pm,n-1}$. 
The $-$ and $+$ superscripts on the left boundary states denote, respectively, values taken from the left and right of the interface.
This distinction is essential, since at cell interfaces the solution of the considered RP can exhibit discontinuities.
Details on how the non-linear system was solved are provided in section 1 of the supplementary material.
%The resulting non-linear system was then analytically manipulated to obtain a non-linear equation that was solved numerically using a standard Newton method. 
The solution to this problem internal to cell $S_i$, namely $\mathbf{Q}_{i-\frac12}^{+,n-1}$, is then used to find the reconstruction polynomial $\mathbf{w}_i^n(x)$.

Due to the structure of our numerical method, the reconstruction polynomials $\mathbf{w}_i^n(x)$ need to be evaluated only in the quadrature nodes, i.e  $x_{i-\frac{1}{2}}$ and $x_{i+\frac{1}{2}}$ for a second-order method, and $x_{i-\frac{1}{2}}, x_i$ and $x_{i+\frac{1}{2}}$ for a third-order method.
In particular, for a second-order method, a first degree polynomial evaluated in the 2 quadrature nodes $x_{i\pm\frac{1}{2}}$ of the computational cell $S_i$ at time $t^n$ reads
\begin{linenomath}
    \begin{equation} \label{eq:rec2}
        \mathbf{w}_i^n(x_{i\pm\frac{1}{2}}) = \mathbf{Q}_i^n \pm \frac{1}{2}(\mathbf{Q}_{i+\frac12}^{-,n-1} - \mathbf{Q}_{i-\frac12}^{+,n-1}).
    \end{equation}
\end{linenomath}
Similarly, for a third-order method, a second degree polynomial evaluated in the 3 quadrature nodes $x_p \in \{x_{i\pm\frac{1}{2}},x_i\}$ of the computational cell $S_i$ at time $t^n$ reads
\begin{linenomath}
    \begin{equation} \label{eq:rec3}
        \mathbf{w}_i^n(x_p) = a+bx_p+cx_p^2,
    \end{equation}
\end{linenomath}
where 
\begin{linenomath}
    \begin{equation}
    a= \mathbf{Q}_{i-\frac12}^{+,n-1},
    \end{equation}
\end{linenomath}
\begin{linenomath}
    \begin{equation}
    b=\frac{2}{\Delta x}(-2\mathbf{Q}_{i-\frac12}^{+,n-1} - \mathbf{Q}_{i+\frac12}^{-,n-1} + 3\mathbf{Q}_i^n),
    \end{equation}
\end{linenomath}
and 
\begin{linenomath}
    \begin{equation}
    c=\frac{3}{(\Delta x)^2}(\mathbf{Q}_{i-\frac12}^{+,n-1} + \mathbf{Q}_{i+\frac12}^{-,n-1} -2\mathbf{Q}_i^n).
    \end{equation}
\end{linenomath}
Both reconstruction polynomials are fully determined by 
either computing a slope or interpolating them through the boundary states $\mathbf{Q}_{i\pm\frac{1}{2}}^{\mp,n-1}$ associated with cell $S_i$, and by enforcing the conservation property, which ensures that the cell average of cell $S_i$ at time $t^n$ is $\mathbf{Q}_i^n$.
Consequently, if the cell average is that of the steady-state solution on the same computational cell, and the boundary states lie on the same steady-state solution, then the resulting reconstruction polynomials evaluated in the quadrature nodes $x_p$ coincide with the steady-state solution.  
We underline that this feature of the GRP-based reconstruction procedure applies only to first and second degree polynomials. Higher degree polynomials require additional information coming from neighboring cells $S_{i\pm1}$ for their complete identification~\cite{GRP}. Thus, they may fail to coincide with the steady-state solution.

\subsection{Summary of the method}
Consider a non-conservative hyperbolic PDE system with source terms written as in~\eqref{eq1}.
The following list summarizes how to compute its numerical solution with the proposed methodology:
\begin{itemize}
    \item[1)] Discretize the spatial domain $\Omega$ into $N$ computational cells $S_i$, and the temporal domain $\mathcal{T}$ into $K$ temporal cells $T^n$.
    \item[2)] At each time step $t^n$ and for each cell $S_i$, compute the steady-state solution of problem~\eqref{eq1} using eq.~\eqref{eq:ssSys1} and eq.~\eqref{eq:ssSys2}.
    \item[3)] At each time step $t^n$ and for each cell $S_i$, compute the spatial reconstruction polynomials using eq.~\eqref{eq:rec2} or eq.~\eqref{eq:rec3}.
    \item[4)] At each time step $t^n$ and for each cell $S_i$, compute the space-time predictions using information from points 2) and 3), and the GRP solver presented in section~\ref{DET}.
    \item[5)] At each time step $t^n$ and for each cell $S_i$, use information from points 2) and 4) to approximate integrals~\eqref{eq:CA},~\eqref{eq:NCs},~\eqref{eq:Sources}, and ~\eqref{eq7}.
    \item[6)] At each time step $t^n$ and for each cell $S_i$, update the solution using eq.~\eqref{eq4}.
    \item [7)] Repeat points 2)-6) till the final simulation time is reached.
\end{itemize}

\section{Numerical tests} \label{Numerical tests}
This section presents the numerical tests conducted to evaluate the performance of the proposed high-order numerical method. 
Two representative test cases are considered: the scalar Burgers’ equation and the hyperbolized BFEs system introduced in~\cite{LTS}. 
To carry out these tests, we employ both second- and third-order implementations of four different methods for comparison.
The first is our method, which combines the GRP-based reconstruction with the well-balanced DET solver (GRP+DET-WB). The second method employs the GRP-based reconstruction, coupled with the original, non-well-balanced DET solver (GRP+DET). The third method combines the WENO reconstruction with the well-balanced DET solver (WENO+DET-WB). The fourth method is a reference approach based on the WENO reconstruction combined with the original DET solver (WENO+DET). More details about the WENO reconstruction can be found in~\cite{WENO}.

\subsection{Scalar case: Burgers' equation}
\begin{figure}[t]
    \centering
    \includegraphics[width=0.99\linewidth]{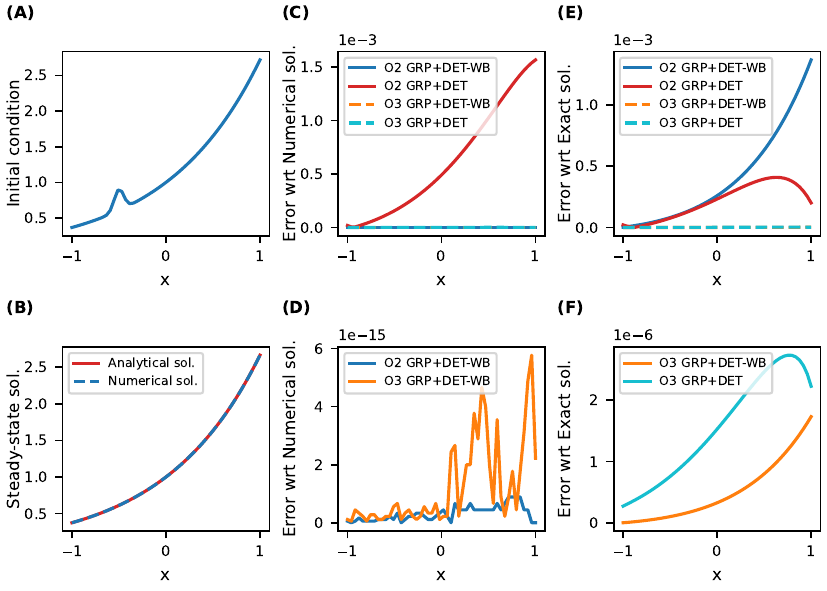}
    \caption{\textbf{Burgers' problem.} 
    %Left column: initial condition (top panel) and exact/numerical reference solutions (bottom panel). Middle column: error between the numerical solution obtained by solving the transient problem (GRP+DET-WB, GRP+DET) and the numerical reference solution. The bottom panel shows a zoom-in on the errors for the second- and third-order GRP+DET-WB schemes. Right column: error between the numerical transient solution and the exact steady-state solution~\eqref{analyticalSol}. The bottom panel shows a zoom-in on the errors for the third-order schemes. 
    Initial condition (A), steady-state solution (B), and errors in space between the numerical solution and the steady-state solution obtained either numerically (C and D) or analytically (E and F). Results are shown for both second- and third-order implementations of the GRP+DET and the GRP+DET-WB methods.
    }
    \label{fig:burgers}
\end{figure}

Let us replicate the test proposed by \citet{GF}, where we consider the following one-dimensional scalar Burgers' equation with algebraic nonlinear source term, expressed in quasi-linear form~\eqref{eq1} as
\begin{linenomath}
\begin{equation} \label{burgers}
    \partial_t q + q\partial_x q = q^2, \quad x\in\Omega, \quad t\in\mathcal{T},
\end{equation}
\end{linenomath}
with initial condition (see figure~\ref{fig:burgers}A) defined as
\begin{linenomath}
\begin{equation}\label{eq:burgIC}
    q\big(x,0\big) = \exp\big(x\big) + 0.3 \exp \big(-200(x+0.5)^2 \big), \quad x\in \Omega, %, C=1.
\end{equation}
\end{linenomath}
and left boundary condition as
\begin{linenomath}
    \begin{equation}
        q(x_A,t) = \exp\big(x_A\big), \quad t\in \mathcal{T}.
    \end{equation}
\end{linenomath}
A transparent right boundary condition is also enforced.
Stationary solutions of eq.~\eqref{burgers} can be computed either analytically as
\begin{linenomath}
\begin{equation} \label{analyticalSol}
    q(x) = \exp \big(x \big), \quad x\in \Omega, %C\in\mathbb{R}   
    \end{equation}
\end{linenomath}
%by solving equation~\eqref{burgers} with the condition $q(x_A) = \exp \big(-x_A \big)$, 
or numerically by applying the procedure described in section~\ref{SteadyState}.

The goal of this test is to demonstrate the ability of our numerical scheme to recover and preserve the stationary solution of eq.~\eqref{burgers} when the initial condition~\eqref{eq:burgIC} is given by a small spatial perturbation of that stationary solution.
% In particular, we replicate the test proposed by \citet{GF}.
We consider a computational domain $\Omega=[-1,1]$, discretized through a mesh of $N=50$ computational cells, and we stop our simulations at the final time $t^K=\SI{40}{s}$. We use a CFL number of 0.9 to compute $\Delta t^n = \text{CFL} \cdot \Delta x / \nu$, with  
\begin{equation}
    \nu = \max(||\mathbf{q}^n||_{L^{\infty}}, ||\mathcal{V}^n||_{L^{\infty}}),
\end{equation}
where $\mathbf{q}$ is the vector containing the cell averages $q_i^n$ of the solution in $S_i$ at time $t^n$, and $\mathcal{V}$ is a vector containing the maximum, on each computational cell, between the left and right states $q_{i-\frac{1}{2}}^{+,n}$ and $q_{i+\frac{1}{2}}^{-,n}$ and their average $\mathcal{M} = 0.5 (q_{i-\frac{1}{2}}^{+,n} + q_{i+\frac{1}{2}}^{-,n})$. The time step $\Delta t^n$ is updated at each iteration of the method and, if necessary, when computing space-time predictions. In this way, we account for the highest wave speed that arises in the presence of shocks and rarefactions. Moreover, we solve equation \eqref{eq18} employing a maximum number of iterations equal to the order of accuracy of the method~\citep{DET}.
% These considerations are not necessary when solving the problem in its conservative form or with a nonlinear reconstruction in Godunov's sense.}
% The initial condition employed for simulations (see figure~\ref{fig:burgers}A) is 
% \begin{linenomath}
% \begin{equation}
%     q = \exp\big(x\big) + 0.3 \exp \big(-200(x+0.5)^2 \big), \ x\in \Omega. %, C=1.
% \end{equation}
% \end{linenomath}

In figures~\ref{fig:burgers}C and \ref{fig:burgers}D, we show the absolute errors between the numerical solution of the steady-state problem
obtained through an appropriate Runge-Kutta scheme, and the numerical solution of eq.~\eqref{burgers} obtained using the GRP+DET-WB and the GRP+DET methods.
Similarly, we display in figures~\ref{fig:burgers}E and \ref{fig:burgers}F the errors between the analytical solution~\eqref{analyticalSol} 
and the numerical solution obtained using the GRP+DET-WB and the GRP+DET methods.

We also perform an empirical convergence test to verify whether our method reaches the expected order of accuracy. We report in table~\ref{tab:burgers} $L^1$ and $L^\infty$ error norms between the numerical solution obtained with the GRP+DET-WB method and the analytical solution~\eqref{analyticalSol}, along with the corresponding empirical convergence rates.

\begin{table}\footnotesize %\small
\caption{$L^1$ and $L^\infty$ error norms, and corresponding empirical convergence rates for a second- and third-order implementation of the numerical scheme. $N$ is the number of computational cells.}
    \centering
    \begin{tabular}{c|cccc|cccc}
    \hline
        N & $L^1$ & $L^\infty$ & $O(L^1)$ & $O(L^\infty)$ & $L^1$ & $L^\infty$ & $O(L^1)$ & $O(L^\infty)$ \\ \hline
        & \multicolumn{1}{l}{Order 2} & & & & \multicolumn{1}{l}{Order 3} \\ \hline
32  & 1.9e-03 & 3.2e-03 & -   & -   & 3.8e-06 & 6.4e-06 & -   & -   \\
64  & 4.9e-04 & 8.4e-04 & 2.0 & 1.9 & 4.8e-07 & 8.3e-07 & 3.0 & 2.9 \\
128 & 1.2e-04 & 2.2e-04 & 2.0 & 2.0 & 6.1e-08 & 1.1e-07 & 3.0 & 3.0 \\
256 & 3.1e-05 & 5.5e-05 & 2.0 & 2.0 & 7.6e-09 & 1.3e-08 & 3.0 & 3.0 \\
512 & 7.8e-06 & 1.4e-05 & 2.0 & 2.0 & 9.6e-10 & 1.7e-09 & 3.0 & 3.0 \\ \hline
    \end{tabular}
    \label{tab:burgers}
\end{table}

\subsection{System case: hyperbolized blood flow equations}
Let us consider the hyperbolized BFEs system proposed in~\cite{LTS}, including now a gravity term. This system of PDEs can be written in quasi-linear form as in eq.~\eqref{eq1}, where we have
\begin{linenomath}
    \begin{equation}
        \mathbf{Q} = [A,q,\psi , A_0, h_0, E_e, E_c, p_r]^T,
    \end{equation}
    \begin{equation}
        \mathbf{A}(\mathbf{Q}) =
        \begin{bmatrix}
        0 & 1 & 0 & 0 & 0 & 0 & 0 & 0 \\
        c^2-u^2 & 2u & \frac{A}{\rho}\partial_{\psi}\zeta & \frac{A}{\rho}\partial_{A_0}\zeta & \frac{A}{\rho}\partial_{h_0}\zeta & \frac{A}{\rho}\partial_{E_e}\zeta & \frac{A}{\rho}\partial_{E_c}\zeta & \frac{A}{\rho} \\
        0 & -1/\varepsilon & 0 & 0 & 0 & 0 & 0 & 0 \\
        0 & 0 & 0 & 0 & 0 & 0 & 0 & 0 \\
        0 & 0 & 0 & 0 & 0 & 0 & 0 & 0 \\
        0 & 0 & 0 & 0 & 0 & 0 & 0 & 0 \\
        0 & 0 & 0 & 0 & 0 & 0 & 0 & 0 \\
        0 & 0 & 0 & 0 & 0 & 0 & 0 & 0 \\
        \end{bmatrix},
    \end{equation}
\end{linenomath}
with $u = q/A$, and $c=\sqrt{\frac{A}{\rho}\partial_{A}\zeta}$, and where
\begin{linenomath}
    \begin{equation}
        \mathbf{S}(\mathbf{Q}) = 
            \big[0, Rq/A+Ag_x, -\psi/\varepsilon, 0, 0, 0, 0, 0\big]^T.
    \end{equation}
\end{linenomath}
$x$ is the axial coordinate along the vessel, $t$ is the time, $A(x,t)$ represents the cross-sectional area of the vessel lumen, $q(x,t)$ is the flow rate, and $\psi(x,t)$ is an auxiliary variable used in the hyperbolization process~\citep{hypBFEs} along with the relaxation time $\varepsilon>0$. $R(<0)$ is the coefficient of the friction term, while $g_x(x)$ is the projection of gravity along the vessel's axis defining the gravity term. $A_0(x)$, $h_0(x)$, $E_e(x)$, $E_c(x)$, and $p_r(x)$ are spatial parameters, assumed to be either constant or continuously varying along the vessel's axis. 
These parameters characterize the vessel's wall mechanics through a nonlinear constitutive relation that links the cross-sectional area to the blood pressure. 
%In this context, pressure and cross-sectional area can thus be used interchangeably, owing to the pressure–area relation that establishes a direct correspondence between them.
%$p(x,t) = p_r(x) + \zeta(x,t)$. 
%These parameters are assumed to be either constant or continuously varying along the vessel's axis.
More details about the spatial parameters, the adopted pressure-area relation, and the system itself can be found in~\citep{LTS}. 
Additionally, details of the eigenstructure of this system are given in~\citep{LTS}, and more extensively in~\citep{hypBFEs, JUNC}.
Finally, we recall that the computation of the boundary states for RP~\eqref{eq:RP11}, associated with the hyperbolized BFEs, is carried out as described in section 1 of the supplementary material.

To test our method on this PDE system, we considered two geometries. First, we used a single blood vessel to verify the empirical convergence of the method to the desired order of accuracy and to compare its efficiency to that of other numerical schemes. Later, we used an arterial network to assess the well-balanced property of the method over a complex geometry and its performance in transient scenarios. Specifically, we considered a reduced version of the anatomically detailed arterial network described in~\cite{ADAN}, ADAN86~\cite{ADAN86}, which is an open network composed of 86 arteries that resemble the real human anatomy.
Among all the arteries of ADAN86, we also considered the right internal carotid artery (ICA) as the single blood vessel on which run the first set of tests, due to its peculiar geometry.

The PDE system was solved in each artery of the ADAN86 network using an initial condition given by constant pressure equal to \SI{60}{mmHg}, and constant zero flow rate.
Additionally, three different types of boundary conditions were enforced. At the inlet of the ascending aorta, an inflow boundary condition was set. At the outlet of terminal vessels, either a Dirichlet boundary condition in cross-sectional area or an outflow boundary condition were used.
Finally, at the joins between two or three arterial segments, junction coupling conditions were prescribed. Further details on the treatment of these boundary conditions and their assignment in the context of the proposed numerical scheme are given in \citep{LTS, JUNC, BC}.
We emphasize that boundary and coupling conditions computed as reported in these references are used by the implemented GRP-based reconstruction, thereby obviating the need for more sophisticated techniques, such as ghost cell filling, to achieve high-order accuracy. 
We also underline that, given the high computational cost of computing~\eqref{eq7} along junctions, we adopt the local time-stepping algorithm proposed in~\citep{LTS}, where, for all tests, the user-defined CFL value is 0.8.

\subsubsection{Efficiency analysis for a single blood vessel (ICA test)}

\begin{figure}[!ht]
\centering
\includegraphics[width = 0.8\textwidth]{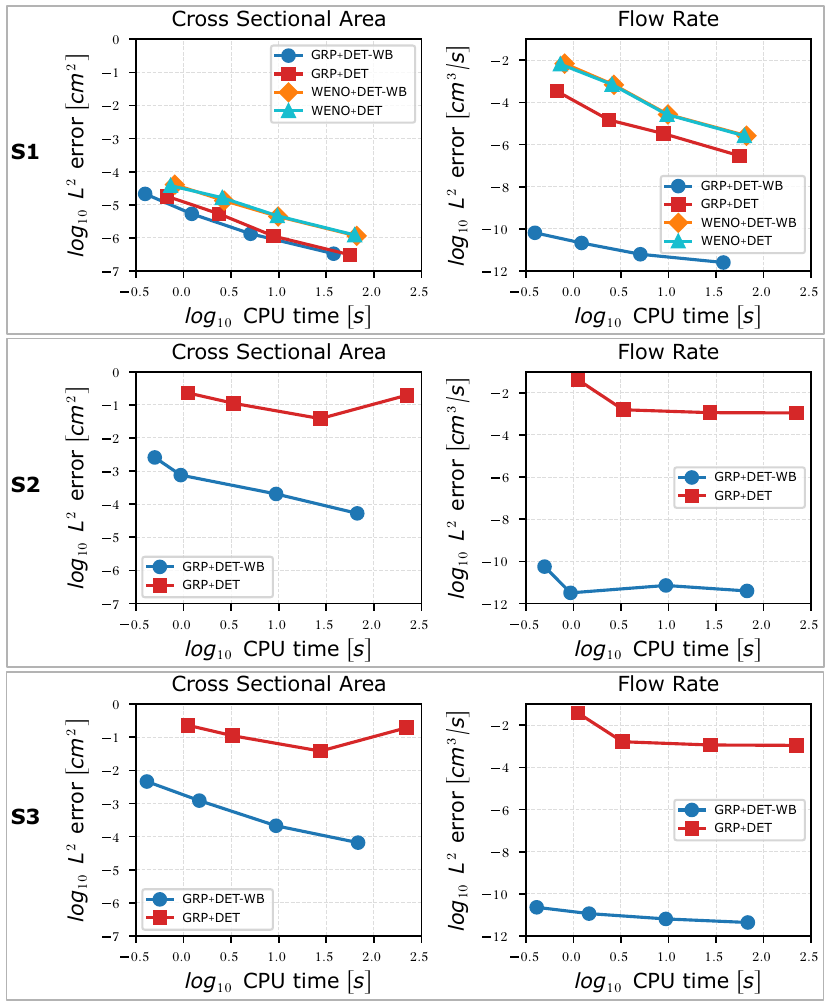}
\caption{\textbf{Efficiency plots for ICA test.} CPU times versus $L^2$ error norms between the numerical solution and either the exact solution (S1 and S2) or a reference solution (S3), for all the considered scenarios and for 4 consecutive mesh refinements. Each row refers to a different scenario (S1, S2, and S3). Results are shown for both cross-sectional area (left panels) and flow rate (right panels) in logarithmic scale on both axes.}
\label{fig:EFF}
\end{figure}

An efficiency analysis was performed to assess if our method was able to attain a prescribed error at low computational cost. 
To this end, we considered three different scenarios of increasing complexity, all defined on the same blood vessel, the ICA, and subject to identical boundary conditions.
The difference among the scenarios arises solely from the choice of both the gravity term and the spatial parameters, in order to highlight the importance of combining appropriate reconstruction techniques with well-balanced solvers to achieve maximal efficiency.

%that highlight the importance of combining appropriate reconstruction techniques with well-balanced numerical schemes to maximize efficiency. 

The first scenario (S1) considers the parameters of the pressure-area relation to be constant along all the vessel axis. 
Specifically, $A_0=0.24$~\si{cm^2}, $h_0=0.05$~\si{cm}, $E_e=3.6\cdot10^6$~\si{Pa}, and $E_c=9\cdot10^8$~\si{Pa}.
It also assumes the gravity projection $g_x$ to be constant and equal to \SI{981}{cm/s^2}. For this first scenario, we compared the results obtained using a second-order implementation of the GRP+DET-WB, the GRP+DET, the WENO+DET-WB, and the WENO+DET methods.

%with our second-order method, which uses the GRP reconstruction combined with the well-balanced version of the DET solver (GRP+DET-WB), to those obtained with the GRP reconstruction combined with the original version of the DET solver (GRP+DET), with the WENO reconstruction combined with the well-balanced version of the DET solver (WENO+DET-WB), and with the WENO reconstruction combined with the original version of the DET solver (WENO+DET). More details about the WENO reconstruction can be found in~\cite{WENO}.

The second scenario (S2) considers the parameters of the pressure-area relation to be variable and continuous in space~\cite{ADAN}. Additionally, a smooth function of the gravity projection is taken into account
\begin{linenomath}
    \begin{equation}
        g_x(x) = |g|\big[\exp(-x) - \exp(-L)\big],
    \end{equation}
\end{linenomath}
where $L$ is the axial length of the ICA, while $|g|=\SI{981}{cm/s^2}$ is the gravitational acceleration modulus.
In this second case, we only compared results obtained with the second-order GRP+DET-WB method, to those obtained with the second-order GRP+DET method, to better highlight the differences among the well-balanced and the non-well-balanced solvers.

The third and final scenario (S3) considers a variation to the second case, where the gravity projection is now a polyline. Here, we took into account the real 3D geometry of the ICA and we projected the gravitational acceleration on the vessel axis. 

For all the scenarios, we assumed a final simulation time of \SI{10}{s}, and a minimum number of 4 spatial computational cells that are doubled for each of the 4 considered mesh refinements. 
Additionally, we 
%assumed constant pressure equal to \SI{60}{mmHg}, and constant zero flow rate as initial conditions. We also 
assumed a no-flow boundary condition at the inlet of the vessel, and fixed cross-sectional area at the outlet such that the corresponding pressure was \SI{60}{mmHg}. In all the cases, we run our tests on a workstation that had a Intel Core i9 processor with 16 cores and 24 threads (\SI{3.2}{GHz} clock speed), and \SI{64}{GiB} of RAM, using 1 thread per test.

Figure~\ref{fig:EFF} shows the obtained results at the final simulation time for a second-order implementation of all the considered methods in terms of CPU time and $L^2$ error norms between the numerical solution and either the exact solution (S1 and S2), or a reference solution (S3) obtained by running a test with 4 times the number of computational cells considered in the final mesh refinement.
For all tests, we show results in terms of cross-sectional area and flow rate.

\subsubsection{Empirical convergence rates for ICA test}

% An empirical convergence test has been made to verify if the expected theoretical order of accuracy of our numerical method is reached. We computed empirical convergence rates for the results of the tests on the ICA obtained with the GRP+DET-WB method for all three considered scenarios (S1, S2, and S3).

% Errors between the computed numerical solutions and either the exact solution (S1 and S2) or a reference solution (S3) have been computed in $L^1$ and $L^{\infty}$ norms. Results are reported in table~\ref{tab1} in terms of cross-sectional area. Results for the flow rate are not given since errors below $10^{-10}$ were always achieved.

An empirical convergence test was carried out to verify whether the proposed method attains its expected theoretical order of accuracy. Convergence rates were computed for the ICA test cases obtained with the GRP+DET-WB method across the three scenarios (S1, S2, and S3). The errors between the computed numerical solutions and either the exact solution (S1 and S2) or a reference solution (S3) were evaluated in both $L^1$ and $L^\infty$ norms. The results, reported in table~\ref{tab1}, are given in terms of cross-sectional area. Errors for the flow rate are omitted, as they consistently remained below $10^{-10}~cm^3/s$.

\begin{table}[ht]\footnotesize
\caption{$L^1$ and $L^\infty$ error norms, and corresponding empirical convergence rates for the cross-sectional area $A(x,t)$ at the final simulation time and for scenarios S1, S2, and S3 for the ICA test, for a second- and a third-order implementation of the numerical scheme. $N$ is the number of computational cells.}
\centering
\begin{tabular}{cc|cccc|cccc}
\hline
Variable & N    & $L^1$  & $L^\infty$ & $O(L^1)$  & $O(L^\infty)$ & $L^1$  & $L^\infty$ & $O(L^1)$  & $O(L^\infty)$\\
\hline
\multicolumn{2}{l|}{Scenario S1} & \multicolumn{1}{l}{Order 2} & & & & \multicolumn{1}{l}{Order 3}\\
\hline
$A$ [$\si{cm^2}$] & 4     & 7.84e-05 & 5.92e-06 & --- & --- & 8.49e-08 & 6.28e-09 & ---  & ---   \\
                  & 8     & 1.96e-05 & 1.48e-06 & 2.0 & 2.0 & 1.08e-08 & 1.17e-09 & 2.97 & 2.42 \\
                  & 16    & 4.90e-06 & 3.72e-07 & 2.0 & 2.0 & 1.36e-09 & 1.68e-10 & 2.99 & 2.80 \\
                  & 32    & 1.23e-06 & 9.30e-08 & 2.0 & 2.0 & 1.74e-10 & 2.34e-11 & 2.97 & 2.84 \\

\hline
\multicolumn{2}{l|}{Scenario S2} & \multicolumn{1}{l}{Order 2} & & & & \multicolumn{1}{l}{Order 3}\\
\hline
$A$ [$\si{cm^2}$] & 4   & 8.59e-03 & 9.59e-04 & --- & --- & 5.90e-03        &  4.37e-04         & ---   & ---   \\
                  & 8   & 2.49e-03 & 2.85e-04 & 1.78 & 1.75 & 9.72e-04 & 1.13e-04     & 2.60  & 1.96  \\
                  & 16  & 6.71e-04 & 7.84e-05 & 1.89 & 1.86 & 1.36e-04  & 2.03e-05   & 2.84  & 2.47  \\
                  & 32  & 1.74e-04 & 2.04e-05 & 1.95 & 1.94 & 1.72e-05  & 2.82e-06    & 2.99  & 2.84  \\

\hline
\multicolumn{2}{l|}{Scenario S3} & \multicolumn{1}{l}{Order 2} & & & & \multicolumn{1}{l}{Order 3}\\
\hline
$A$ [$\si{cm^2}$] & 4   & 1.64e-02 & 1.64e-03 & --- & --- & 9.79e-04        & 1.16e-04         & ---   & ---   \\
                  & 8   & 4.04e-03 & 5.15e-04 & 2.02 & 1.67 & 8.38e-04        & 8.36e-05         & 0.22  & 0.47  \\
                  & 16  & 6.66e-04 & 9.75e-05 & 2.60 & 2.40 & 7.48e-05        & 1.36e-05         & 3.49  & 2.62  \\
                  & 32  & 2.05e-04 & 3.70e-05 & 1.70 & 1.40 & 3.17e-05        & 4.34e-06         & 1.24  & 1.65  \\
\hline
    \end{tabular}
    \label{tab1}
\end{table}

\subsubsection{Well-balance property for the ADAN86 geometry (deadman test)}
This test is designed to verify if our method can accurately approximate a specific steady-state solution, characterized by zero flow throughout the network and a hydrostatic pressure distribution (i.e. the pressure varies linearly along the vertical axis of the body due to gravity, with higher values in lower regions, such as in the legs, and lower values in elevated regions, such as in the head).
% The well-balanced property of a numerical method is known to be important for preserving steady state solutions. In the case of BFEs, requiring a method to be well-balanced is equivalent to asking the method to be able to correctly retrieve the solution known as deadman. 
%The deadman solution is a trivial solution given by a zero flow condition everywhere along the network and an hydrostatic pressure distribution.

% In order to verify that our method is able to correctly approximate the deadman solution, we decided to perform the following test. 
To this end, we considered the ADAN86 network in the upright posture, enabling the simulation of hydrostatic pressure gradients induced by gravity. Then, we assumed the parameters of the pressure-area relation to be variable and continuous in space. 
In this test, the gravity term $g_x(x)$ was a polyline. Specifically, we took into account the real 3D geometry of the different vessels and we projected the gravitational acceleration on the vessels' axes.
%the projection of the gravitational acceleration along the vessels' axes. 
All blood vessels were discretized using a maximum mesh spacing of \SI{1}{cm}, and the simulations were performed over \SI{20}{s}. 
%We further assumed that the gravity projection $g_x(x)$ corresponds to the true component of gravitational acceleration along the vessel axes, under the condition that the individual is in a static, upright posture.
% We also assumed the gravity projection $g_x(x)$ to be the actual projection of the gravitational acceleration along the vessels axes, as if the body of the person we are considering is standing still. 
%We then assumed to discretize all the blood vessels with a maximum mesh spacing of \SI{1}{cm}, and to run a simulation for \SI{20}{s} with CFL=0.8. 
%We considered an initial condition given by no-flow and constant pressure (\SI{60}{mmHg}) everywhere in all the vessels of the network. 
Finally, we assumed no-flow boundary conditions at both inlet and outlets of the terminal vessels.
%while we solved junction coupling conditions as in~\cite{LTS, JUNC}.

Results are shown in figure~\ref{fig:deadman} and figure~2 of the supplementary material in terms of flow rate (top panel) and pressure (bottom panel). Specifically, the two networks in each panel display the errors at the final simulation time along the ADAN86 network between the deadman solution and the numerical solutions obtained with a second (and third) order implementation of the GRP+DET-WB and the GRP+DET methods. A focus on three vessels (left anterior cerebral artery, thoracic aorta, and left femoral artery) is provided in the middle of each panel, where results are reported for the midpoint of the vessels axis and for the last second of simulation.
Additionally, figure~\ref{fig:hydro} and figure~3 of the supplementary material show the obtained pressure distribution along the network for both the applied methods of order two and three at the final simulation time. A focus on different vessels at different heights of the network is provided in the middle of the panel, where results are displayed for the midpoint of the chosen vessels axis.

\begin{figure}[p]
\centering
\includegraphics[width=0.775\textwidth]{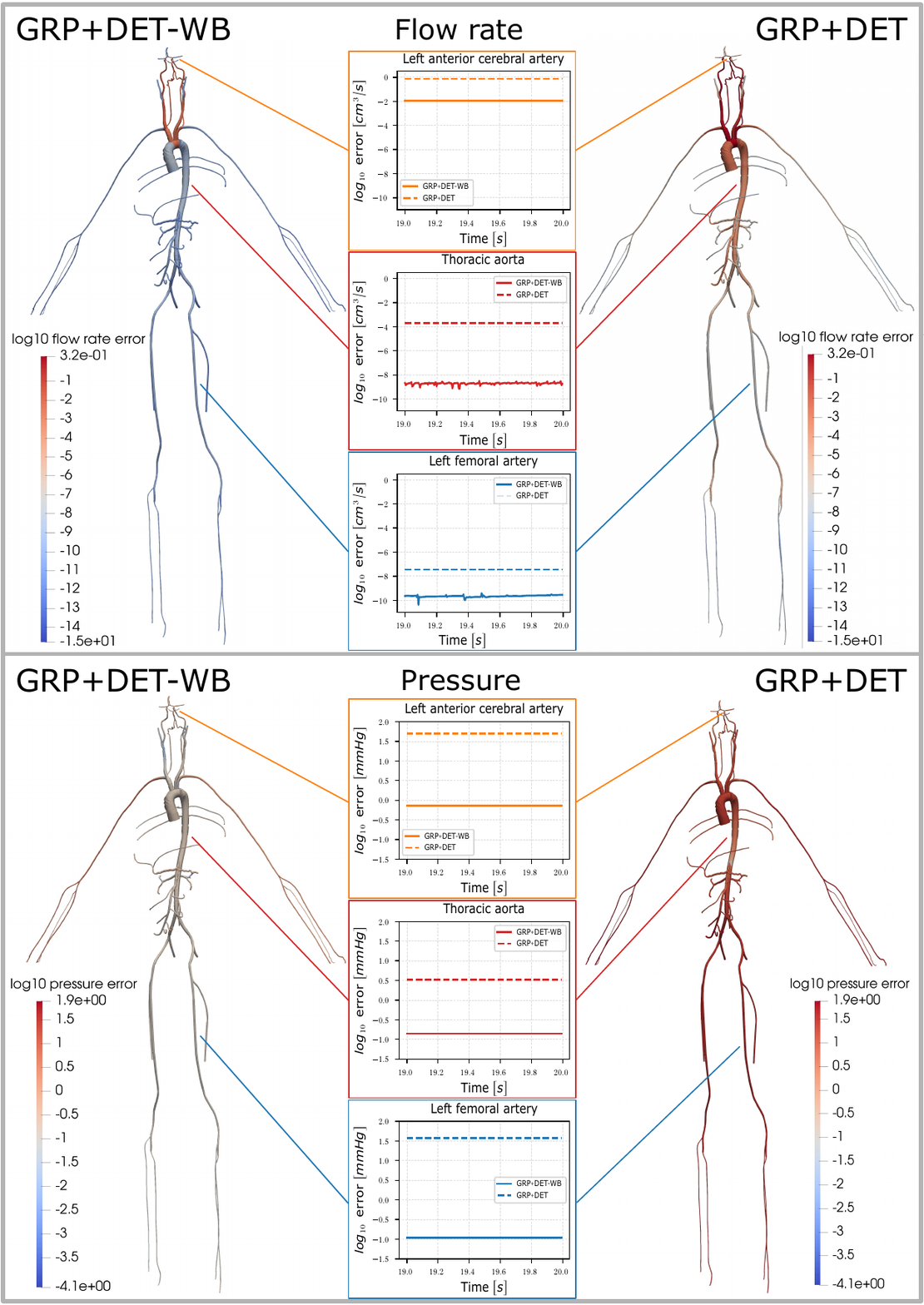}
\caption{\textbf{Deadman test.} A representation of the ADAN86 network is shown. In the top panel, the colors indicate the errors between zero-flow solution and the second-order numerical solution computed with either the GRP+DET-WB, or the GRP+DET. In the bottom panel, the colors indicate the errors between the reference hydrostatic pressure distribution and the numerical solution computed with either the GRP+DET-WB, or the GRP+DET. A focus on three vessels is provided in the middle of both panels. All the results are shown in logarithmic scale.}
\label{fig:deadman}
\end{figure}

\begin{figure}[!ht]
\centering
\includegraphics[width=0.99\textwidth]{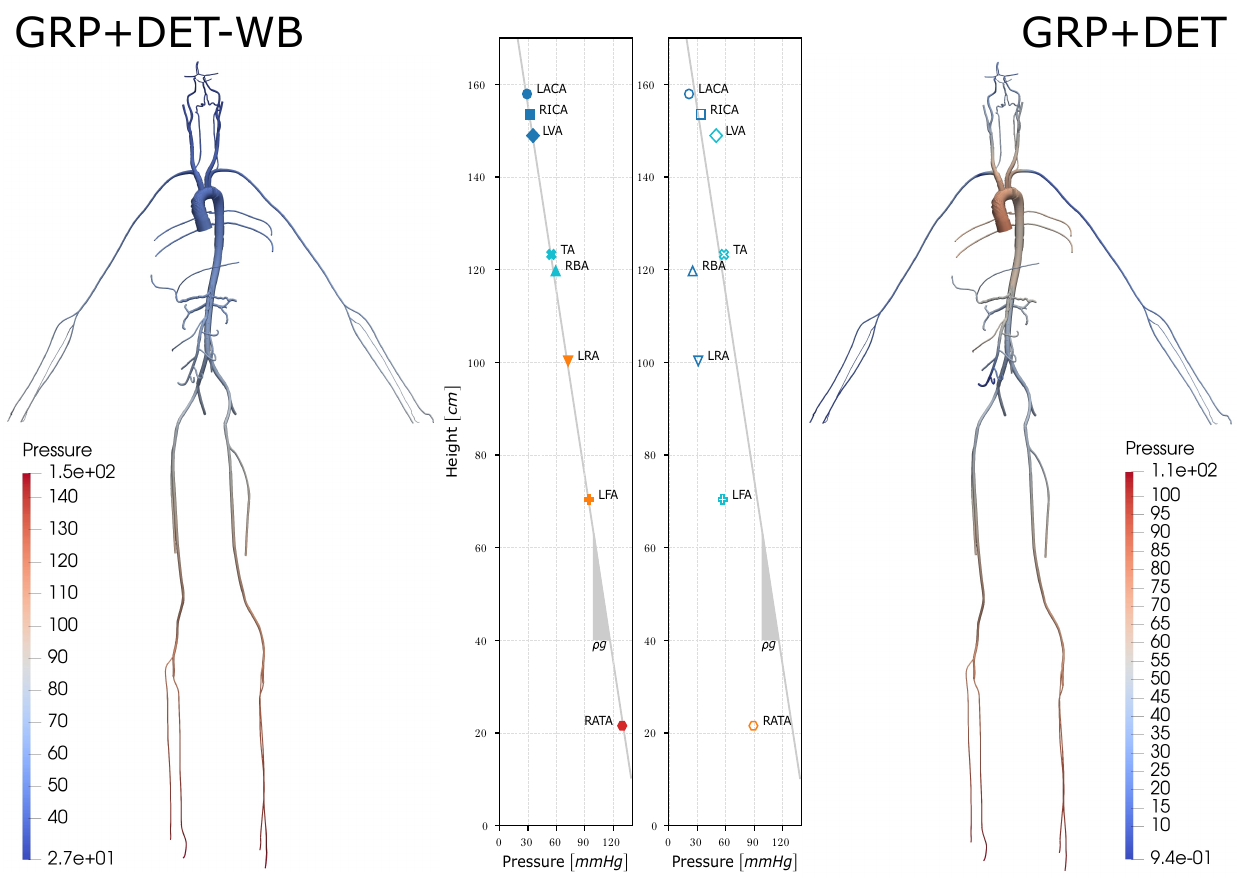}
\caption{\textbf{Pressure distribution.} Pressure distribution at the final simulation time along the ADAN86 network obtained with a second-order implementation of both the GRP+DET-WB method (left) and the GRP+DET method (right). A focus on eight vessels is provided in the middle of the panel, showing how the GRP+DET-WB results respect the expected hydrostatic distribution indicated in light gray. The considered vessels are: LACA: left anterior cerebral artery, RICA: right internal carotid artery, LVA: left vertebral artery, TA: thoracic aorta, RBA: right brachial artery, LRA: left radial artery, LFA: left femoral artery, RATA: right anterior tibial artery.}
\label{fig:hydro}
\end{figure}

\subsubsection{Transient simulation for ADAN86}
The final test we present aims to assess the capability of our method to accurately capture transient solutions in complex geometries. 
%when applied to complex geometries to describe transient solutions. 
Specifically, we wanted to test its ability in retrieving pressure and flow rate curves that are commonly observed when modeling the cardiovascular system.

We considered the ADAN86 network in the upright posture and we assumed the parameters of the pressure-area relation to be variable and continuous in space. As in the previous test, the gravity projection $g_x(x)$ was defined as the projection of the gravitational acceleration along the vessels' axes. All blood vessels were discretized using a maximum mesh spacing of \SI{2}{cm}, and the simulations were performed over \SI{10}{s}.
A periodic inflow boundary condition was prescribed at the aortic inlet~\cite{ADAN} to simulate the behavior of the heart (here not present), thereby reproducing the dynamics of the human cardiac cycle.
Finally, an outflow boundary condition was enforced by coupling terminal vessels to lumped parameters models that mimic the resistance and compliance of the venous system~\cite{LTS}.
We underline that the pressure of the venous system was assumed to be hydrostatically distributed with respect to the vertical axis of the body. In particular, we considered the aortic root to be our reference level.
%We considered a no-flow and constant pressure (\SI{60}{mmHg}) initial condition, and we set an inflow boundary condition at the inlet of the aorta~\cite{ADAN}. We also solved junction boundary conditions as in~\cite{JUNC}, and we connected the outlets of terminal vessels to RCR circuits solved as in~\cite{LTS}.

Results are shown in figure~\ref{fig:transient} and figure~4 of the supplementary material in terms of flow rate and pressure for three specific vessels (left anterior cerebral artery, thoracic aorta, and left femoral artery) at their midpoint along the final cardiac cycle of simulation. 
A reference solution, obtained by running the test with gravity set to zero, is also displayed.

\begin{figure}[t]
\centering
\includegraphics[width=0.99\textwidth]{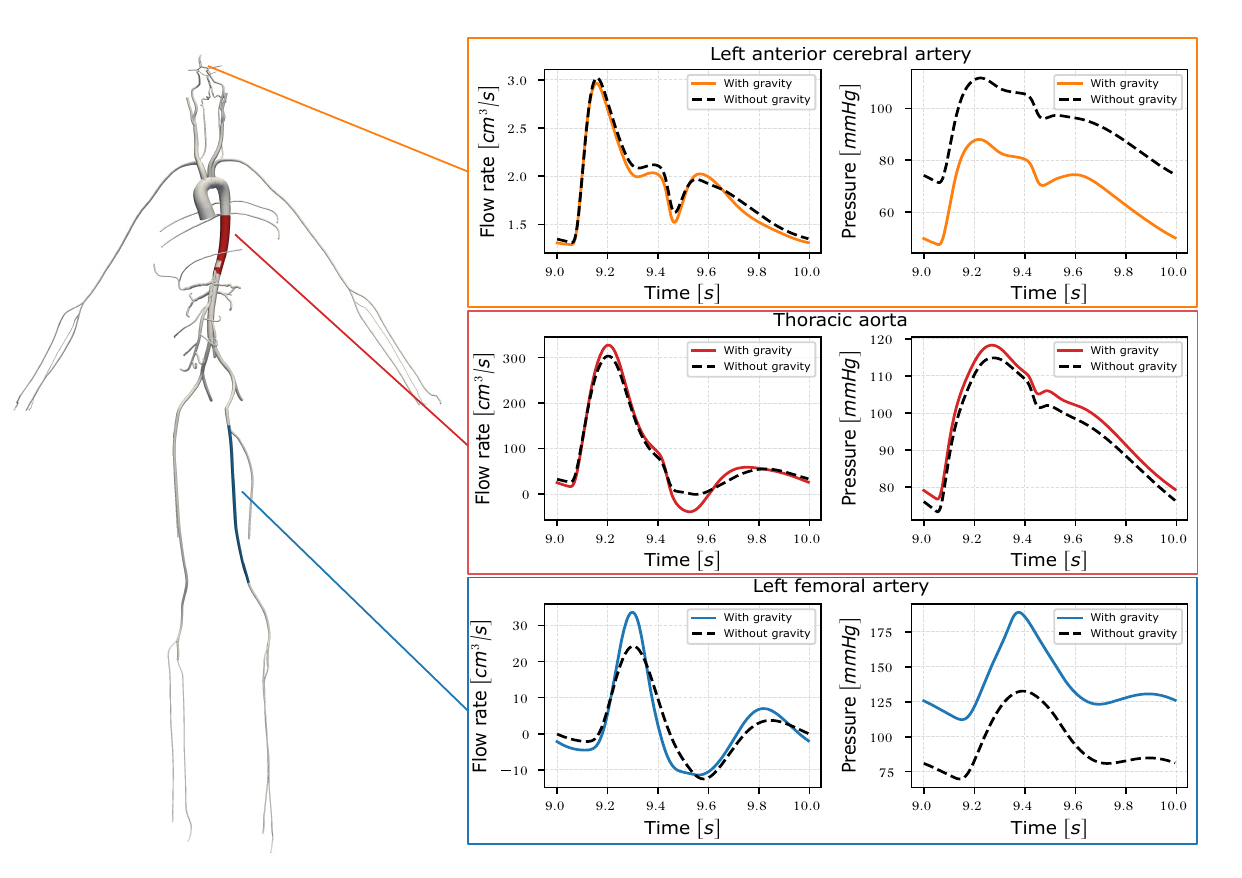}
\caption{\textbf{Transient test.} Flow rate (left column) and pressure (right column) curves along the final cardiac cycle of simulation at the midpoint of three selected blood vessels: left anterior cerebral artery (orange), thoracic aorta (red), and left femoral artery (blue).
Indication on the location of the vessels is shown in the ADAN86 network representation on the left. Black dashed lines represent the reference solution with no gravity in the different vessels.}
\label{fig:transient}
\end{figure}

\section{Discussion} \label{results}
In this section, we discuss results obtained through the numerical tests outlined in section \ref{Numerical tests}, assessing the accuracy and efficiency of our numerical method, and evidencing its limitations.

Empirical convergence results reported in tables \ref{tab:burgers} and \ref{tab1} demonstrate that the GRP+DET-WB scheme converges with the expected order of accuracy to either the analytical steady-state solution or to a reference solution computed on a sufficiently fine mesh.
Additionally, results reported in figure \ref{fig:burgers}D confirm that the method is well-balanced when applied to find the steady-state solution of problem~\eqref{burgers}.
This proves numerically that, as opposed to the standard WENO reconstruction, the GRP-based reconstruction is well-balanced up to the third order of accuracy by construction. 
Indeed, as mentioned in section \ref{rec}, reconstruction polynomials are computed by enforcing that the boundary states associated with each cell are used either to compute the polynomials slopes or as interpolation points, and that the conservation property holds. Consequently, if the cell average is that of the steady-state solution on the same computational cell, and the boundary states lie on
the same steady-state solution, then they are not perturbed by the reconstruction procedure.

This is also confirmed from the efficiency plots reported in figure \ref{fig:EFF}. When considering scenario S1 (top row), we can notice that the choice of the spatial reconstruction plays a significant role: results obtained through the GRP+DET method %\cite{GRP} 
present lower errors, with the same spatial discretization, than results obtained through the WENO+DET method. 
Moreover, the former are comparable, when considering the cross-sectional area, to results obtained through the GRP+DET-WB scheme. As a consequence, we can conclude that the only well-balancing errors introduced within the GRP+DET framework are those associated with the DET solver.

A comparison between numerical results and either the analytical solution of the problem or the reference solution computed on a sufficiently fine mesh (see figures \ref{fig:burgers}E, \ref{fig:burgers}F, and \ref{fig:EFF}), shows that, while the GRP+DET-WB and GRP+DET methods perform similarly in simple cases such as the Burgers' equation and scenario S1 for the BFEs, there is an evident advantage in choosing a well-balanced scheme for more complex scenarios. When considering scenarios S2 and S3 for the BFEs, 
numerical errors for the cross-sectional area and the flow rate obtained through the GRP+DET-WB scheme are significantly lower than those obtained in the GRP+DET setup with analogous meshes, highlighting the higher efficiency of the well-balanced setup. 
% As for the flow-rate, this difference is even more marked, reaching 8 orders of magnitude.
The similar performance of the GRP+DET-WB and GRP+DET schemes in simple scenarios suggests that the approximation of steady-state solutions (see section \ref{SteadyState}) has a significant impact on the overall numerical errors introduced by the scheme.
This poses the problem of choosing an ordinary differential equation solver of high enough accuracy to avoid the introduction of numerical errors that are comparable or higher than those introduced by the absence of well-balancing.

The ability of our numerical method to preserve stationary solutions is also particularly evident when the method is applied to solve the BFEs over a complex geometry like the ADAN86 network. In fact, observing figure~\ref{fig:deadman}, top panel, it is evident that the errors in flow rate between the exact solution given by zero-flow and the numerical solution are always higher for the GRP+DET case with respect to the GRP+DET-WB case. A minimum of two orders of magnitude difference between the solutions obtained with the two methods is always detected. Additionally, we observe that the maximum errors for the GRP+DET-WB case are found in the neck and head region, due to its intricate geometry. This region presents multiple vessels junctions, and blood vessels with high geometrical variability, which results in a gravity projection ranging from \SI{-981}{cm/s^2} to \SI{+981}{cm/s^2} along the same vessel axis that can significantly influence the approximation of the solution.
Similarly, observing the bottom panel of figure~\ref{fig:deadman}, we note that the
%A similar behavior is also present when observing the bottom panel of figure~\ref{fig:deadman}. The
errors in pressure between the deadman solution and the numerical solution are generally larger for the GRP+DET case with respect to the GRP+DET-WB case. These errors strongly affect the representation of the pressure distribution. In particular, we expect to see an hydrostatic pressure distribution along the network when a steady-state is reached, with a maximum pressure in the legs region that gradually decreases towards the cerebral region. 
This behavior is correctly reproduced in figure~\ref{fig:hydro} by the GRP+DET-WB scheme, whereas the GRP-DET scheme fails to capture the hydrostatic distribution, leading to an incorrect approximation of the target solution with a maximal error of \SI{100}{mmHg}.
%While the hydrostatic distribution is visible in figure~\ref{fig:hydro} for the GRP+DET-WB case, it does not emerge for the GRP-DET case, indicating a wrong approximation of the sought solution, with a maximal error of \SI{100}{mmHg}.

%Our method has demonstrated its effectiveness in retrieving and preserving stationary solutions on the ADAN86 network. 
Finally, from the pressure and flow-rate waveforms reported in figure~\ref{fig:transient}, we observe that the effects of the gravity term are most apparent in the pressure curves. The flow-rate waveforms obtained with gravity are in good agreement with those computed without gravity, with only slight differences in curvature and peak values. In contrast, the pressure waveforms with gravity exhibit a vertical shift relative to those without gravity. This shift reflects the position of the vessel under consideration with respect to the aortic root, which serves as the reference level for the assumed hydrostatic pressure distribution of the venous system. These findings are consistent with previous observations reported in~\citep{Fois, MurilloG}, confirming that the proposed method is capable of performing transient simulations in complex scenarios, allowing an accurate reproduction of the behavior of the cardiovascular system.

%it is clear that it is also able to run transient simulations, allowing an accurate reproduction of the behavior of the cardiovascular system. 

%the good agreement between pressure and flow rate waveforms reported in figure~\ref{fig:transient} and the considered reference solution demonstrates that the method can also correctly capture transient solutions, allowing an accurate reproduction of the behavior of the cardiovascular system. 

\section{Conclusion}\label{Concl}

In this work, we proposed a high-order well-balanced numerical scheme for the solution of non-conservative systems of hyperbolic PDEs with source terms. 
The method is based on the GRP reconstruction by~\citet{GRP}, which is well-balanced by construction up to third order of accuracy, provided that boundary states belong to the same family of the steady-state solution of the considered PDEs and the conservation property holds.
In addition, it incorporates a modification of the DET solver~\cite{DET} along the same lines as those proposed by~\citet{GF}. 
Numerical tests performed on the Burgers' equation and on the hyperbolized BFEs~\cite{LTS} with friction, gravity and variable geometrical properties demonstrate the well-balanced property of the scheme and highlight its accuracy and efficiency in complex scenarios. 

However, it is important to mention that in this work we employ the GRP-based reconstruction without making use of any limiter, which is linear in Godunov's sense. As a consequence, spurious oscillations might arise in the presence of discontinuities. A possible strategy to overcome this issue would be to employ an a-posteriori limiter such as the Multi-dimensional Optimal Order Detection (MOOD) strategy~\cite{MOOD}.

\section*{Acknowledgments}
C.C. and L.O.M acknowledge funding by the European Union under NextGenerationEU, Mission 4, Component 2 - PRIN 2022 (D.D. 104/22), project title: Immersed methods for multIscale and multiphysics problems, CUP: E53D2300592 0006. \\
C.C. and A.S. acknowledge the Italian Ministry of Education, Universities and Research (MUR), in the framework of the project DICAM-EXC (Departments of Excellence 2023-2027, grant L232/2016). \\
C.C., C.D., and L.O.M. are members of the ”Gruppo Nazionale per il Calcolo Scientifico" dell’ Istituto Nazionale di Alta Matematica (INdAM-GNCS, Italy).

% \clearpage
% \section{Conclusions}

\bibliographystyle{elsarticle-num-names} %-harv tiene elenco referenze in ordine alfabetico; num-names in ordine di apparizione
\bibliography{main.bib}

\section*{Supplementary material}
\subsection{Solution to RP (23) for the BFEs}
Let us consider the following classical RP for the hyperbolized BFEs presented in the main manuscript at eq. (34)-(36)
\begin{equation}
\partial_t\mathbf{Q}+\mathbf{A}(\mathbf{Q})\partial_x\mathbf{Q}=\mathbf{0}, \quad x\in\mathbb{R}, t>0
\end{equation}
with initial condition
\begin{equation}
    \mathbf{Q}(x,0) = 
    \begin{cases}
        \mathbf{Q}_L, & \text{if }x<x_{i+\frac12},\\
        \mathbf{Q}_R, & \text{if }x>x_{i+\frac12}.
    \end{cases}
\end{equation}
$\mathbf{Q}_L$ and $\mathbf{Q}_R$ are constant states defined to both side of the cell interface $x_{i+\frac12}$.

Our goal is to solve this RP under the condition that the flow is subcritical.
Assuming a subcritical flow is equivalent to asking that the fluid velocity $u$ is smaller than the wave propagation speed $c$, namely $\lambda_1<0$ and $\lambda_8>0$.
Indeed, we recall that the hyperbolized BFEs present 8 waves~\citep{hypBFEs}, with associated eigenvalues given by
\begin{equation}
    \lambda_1(\mathbf{Q})=u-c, \quad \lambda_{2,\dots,7}(\mathbf{Q})=0, \quad \lambda_8(\mathbf{Q})=u+c. 
\end{equation}

\begin{figure}[t]
\centering
\includegraphics[width=0.775\textwidth]{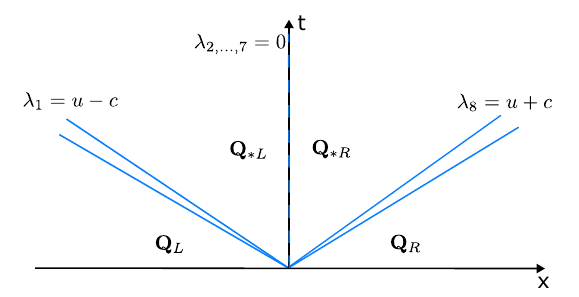}
\caption{\textbf{Wave configuration.} Representation of the wave configuration in the $x-t$ half plane.}
\label{fig:1SM}
\end{figure}

Within the subcritical regime, there is only a possible wave configuration in the $x-t$ half plane (see figure~\ref{fig:1SM}) given by 4 constant states, namely $\mathbf{Q}_L, \mathbf{Q}_{*L}, \mathbf{Q}_{*R},\mathbf{Q}_R$~\citep{spilimbergo2021one}.
$\mathbf{Q}_L$ is left constant state (known) that is separated from the state $\mathbf{Q}_{*L}$ by the left family of waves associated to $\lambda_1$. Similarly, $\mathbf{Q}_R$ is right constant state (known) that is separated from the state $\mathbf{Q}_{*R}$ by the right family of waves associated to $\lambda_8$. Both left and right families of waves can generate either shock waves or rarefactions~\citep{hypBFEs}. Finally, $\mathbf{Q}_{*L}$ and $\mathbf{Q}_{*R}$ are the two unknown states that are separated by contact discontinuities. In order to identify these two unknown states, and thus solve the RP, we need to establish appropriate jump conditions across the family of waves to connect the unknown states $\mathbf{Q}_{*L}$, $\mathbf{Q}_{*R}$ between them and to the initial conditions $\mathbf{Q}_{L}$, and $\mathbf{Q}_{R}$.

The solution to this RP can be obtained using various methods~\citep{Toro}. 
Here we apply a two-rarefaction approximate Riemann solver. Specifically, we assume that the waves associated to $\lambda_1$ and $\lambda_8$ are rarefactions, and we use the related Riemann invariants to identify the unknowns~\citep{hypBFEs, spilimbergo2021one}.
The application of the solver to our problem yields
\begin{equation}\label{eq4SM}
    A_{0,L}=A_{0,*L}, \quad h_{0,L}=h_{0,*L}, \quad E_{e,L}=E_{e,*L}, \quad E_{c,L}=E_{c,*L}, \quad p_{r,L}=p_{r,*L},
\end{equation}
\begin{equation}\label{eq5SM}
    A_{0,R}=A_{0,*R}, \quad h_{0,R}=h_{0,*R}, \quad E_{e,R}=E_{e,*R}, \quad E_{c,R}=E_{c,*R}, \quad p_{r,R}=p_{r,*R},
\end{equation}
and
\begin{equation}\label{eq6SM}
    \begin{cases}
        \psi_L + A_L/\varepsilon = \psi_{*L}+A_{*L}/\varepsilon, \\
        u_{*L} = u_L - \int_{A_L}^{A_{*L}} \frac{\tilde{c}_T(\xi)}{\xi}d\xi, \\
        \psi_R + A_R/\varepsilon = \psi_{*R}+A_{*R}/\varepsilon, \\
        u_{*R} = u_R + \int_{A_R}^{A_{*R}} \frac{\tilde{c}_T(\xi)}{\xi}d\xi, \\
        q_{*L} = q_{*R}, \\
        p_{*L} + \frac12\rho u_{*L}^2 = p_{*R} + \frac12\rho u_{*R}^2.    
    \end{cases}
\end{equation}
The equivalences in eq.~\eqref{eq4SM} and \eqref{eq5SM} identify the first set of unknowns, while the system in eq.~\eqref{eq6SM} leads us to the determination of the remaining six variables, namely $A_{*L/R}$, $q_{*L/R}$, and $\psi_{*L/R}$.

In order to solve system~\eqref{eq6SM}, we first apply the midpoint rule to the two integrals appearing in it.
Specifically, we have
\begin{equation}
    \int_{\hat{A}}^{A} \frac{\tilde{c}_T(\xi)}{\xi}d\xi \cong (A-\hat{A}) \Bigg[\frac{\tilde{c}_T\big(\frac{A+\hat{A}}{2}\big)}{\frac{A+\hat{A}}{2}}  \Bigg], 
\end{equation}
with $\hat{A}=A_L,A_L$, and $A=A_{*L},A_{*R}$.
Later, we assume that the integrand function between $\hat{A}$ and $A$ is flat enough to write that
\begin{equation}
    \frac{\tilde{c}_T\big(\frac{A+\hat{A}}{2}\big)}{\frac{A+\hat{A}}{2}} \approx \frac{\tilde{c}_T(\hat{A})}{\hat{A}}.
\end{equation}
As a result, system~\eqref{eq6SM} becomes
\begin{equation}\label{eq9SM}
    \begin{cases}
        \psi_L + A_L/\varepsilon = \psi_{*L}+A_{*L}/\varepsilon, \\
        u_{*L} = u_L + \tilde{c}_T(A_L) - (A_{*L}/A_L)\cdot \tilde{c}_T(A_L), \\
        \psi_R + A_R/\varepsilon = \psi_{*R}+A_{*R}/\varepsilon, \\
        u_{*R} = u_R - \tilde{c}_T(A_R) + (A_{*R}/A_R) \cdot\tilde{c}_T(A_R), \\
        q_{*L} = q_{*R}, \\
        p_{*L} + \frac12\rho u_{*L}^2 = p_{*R} + \frac12\rho u_{*R}^2.    
    \end{cases}
\end{equation}
We observe that the first 4 equations of the system can all be written in terms of $A_{*L}$ and $A_{*R}$, namely 
\begin{equation}\label{eq10SM}
    \begin{cases}
    \psi_{*L} = g_1(A_{*L}), \\
    \psi_{*R} = g_2(A_{*R}), \\
    u_{*L} = h_1(A_{*L}), \\
    u_{*R} = h_2(A_{*R}), \\
        q_{*L} = q_{*R}, \\
        p_{*L} + \frac12\rho u_{*L}^2 = p_{*R} + \frac12\rho u_{*R}^2,    
    \end{cases}
\end{equation}
where $g_1,g_2,h_1,h_2$ are known.
Additionally, we also note that $q=Au$, which allows us to write the fifth equation as
\begin{equation}
    A_{*L}u_{*L} = A_{*R}u_{*R},
\end{equation}
and to identify a relation between $A_{*L}$ and $A_{*R}$.
Specifically, the previous equation represents an ellipse of the form
\begin{equation}
\frac{\tilde{c}_T(A_R)}{A_R}A_{*R}^2 + (u_R-\tilde{c}_T(A_R))A_{*R}
    +\frac{\tilde{c}_T(A_L)}{A_L}A_{*L}^2 - (u_L+\tilde{c}_T(A_L))A_{*L}=0.
\end{equation}
Writing thus $A_{*R}$ as a function of $A_{*L}$, i.e. $A_{*R}=f(A_{*L})$, which represent the solution to the ellipse, we can now write all the equations of system~\eqref{eq10SM} in terms of the same unknown
\begin{equation}\label{eq13SM}
    \begin{cases}
    \psi_{*L} = g_1(A_{*L}), \\
    \psi_{*R} = g_2(f(A_{*L})), \\
    u_{*L} = h_1(A_{*L}), \\
    u_{*R} = h_2(f(A_{*L})), \\
        A_{*R}=f(A_{*L}), \\
        p_{*L} + \frac12\rho h_1^2(A_{*L}) = p_{*R} + \frac12\rho h_2^2(f(A_{*L})).    
    \end{cases}
\end{equation}
Finally, after this long manipulation of the original system~\eqref{eq6SM}, we observe that system~\eqref{eq13SM} is characterized by a single non-linear equation (the last one) in terms of the unknown $A_{*L}$, and multiple equivalences.
Thus, applying a standard Newton method to the last equation, it allows us to identify $A_{*L}$ and successively to retrieve the remaining unknowns.

\subsection{Additional results related to the third-order implementation of the method}

\begin{figure}[!ht]
\centering
\includegraphics[width=0.775\textwidth]{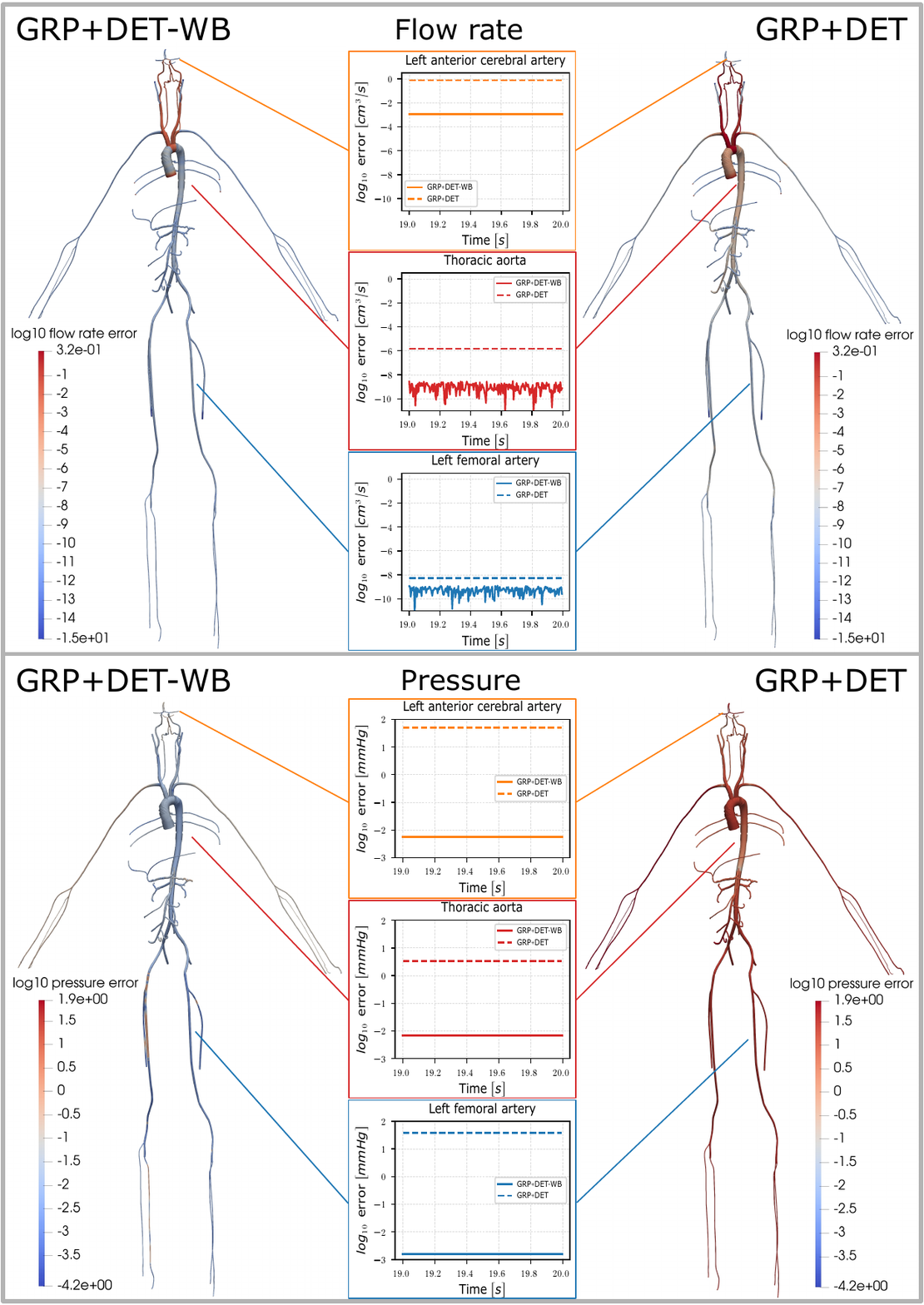}
\caption{\textbf{Deadman test.} A representation of the ADAN86 network is shown. In the top panel, the colors indicate the errors between zero-flow solution and the third-order numerical solution computed with either the GRP+DET-WB, or the GRP+DET. In the bottom panel, the colors indicate the errors between the reference hydrostatic pressure distribution and the numerical solution computed with either the GRP+DET-WB, or the GRP+DET. A focus on three vessels is provided in the middle of both panels. All the results are shown in logarithmic scale.}
\label{fig:deadman3}
\end{figure}

\begin{figure}[!ht]
\centering
\includegraphics[width=0.99\textwidth]{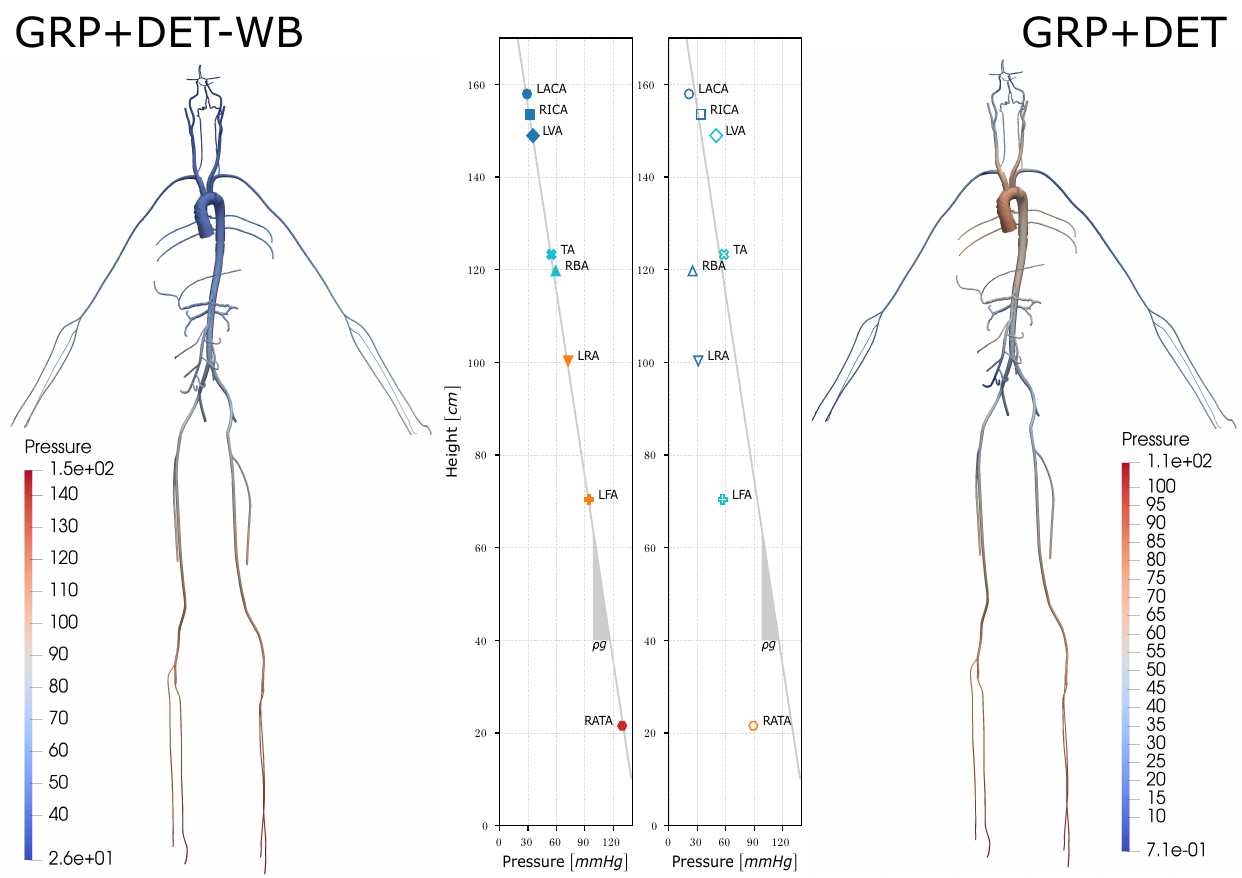}
\caption{\textbf{Pressure distribution.} Pressure distribution at the final simulation time along the ADAN86 network obtained with a third-order implementation of both the GRP+DET-WB method (left) and the GRP+DET method (right). A focus on eight vessels is provided in the middle of the panel, showing how the GRP+DET-WB results respect the expected hydrostatic distribution indicated in light gray. The considered vessels are: LACA: left anterior cerebral artery, RICA: right internal carotid artery, LVA: left vertebral artery, TA: thoracic aorta, RBA: right brachial artery, LRA: left radial artery, LFA: left femoral artery, RATA: right anterior tibial artery.}
\label{fig:hydro3}
\end{figure}

\begin{figure}[!ht]
\centering
\includegraphics[width=0.99\textwidth]{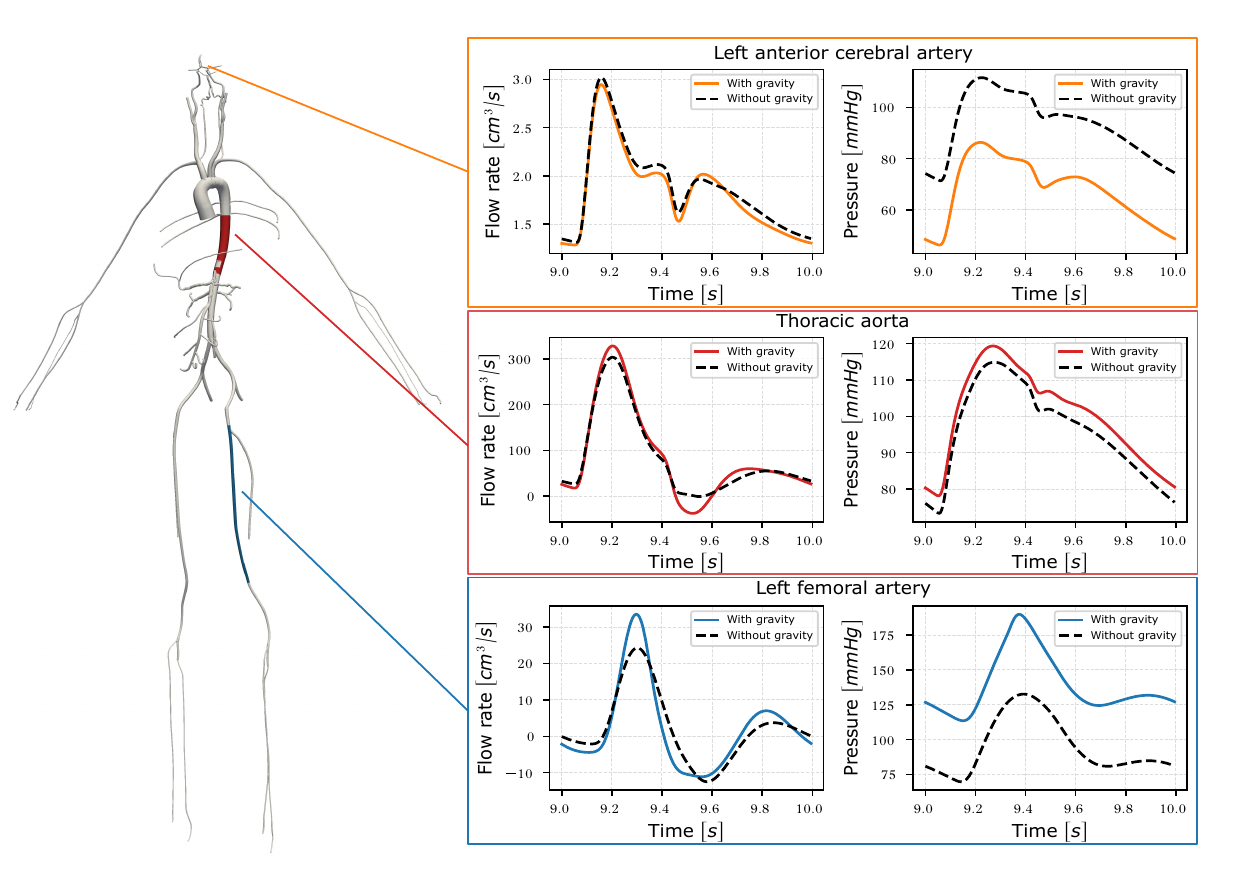}
\caption{\textbf{Transient test.} Flow rate (left column) and pressure (right column) curves along the final cardiac cycle of simulation at the midpoint of three selected blood vessels: left anterior cerebral artery (orange), thoracic aorta (red), and left femoral artery (blue).
Indication on the location of the vessels is shown in the ADAN86 network representation on the left. Black dashed lines represent the reference solution with no gravity in the different vessels.}
\label{fig:transient3}
\end{figure}

% %% Use \subsubsection, \paragraph, \subparagraph commands to 
% %% start 3rd, 4th and 5th level sections.
% %% Refer following link for more details.
% %% https://en.wikibooks.org/wiki/LaTeX/Document_Structure#Sectioning_commands

% %% The Appendices part is started with the command \appendix;
% %% appendix sections are then done as normal sections
% \appendix
% \section{Example Appendix Section}


% %% For citations use: 
% %%       \citet{<label>} ==> Lamport [21]
% %%       \citep{<label>} ==> [21]
% %%
% Example citation, See \citet{lamport94}.

% %% If you have bib database file and want bibtex to generate the
% %% bibitems, please use
% %% \bibliographystyle{elsarticle-num-names} 
% %% \bibliography{biblio}

% %% else use the following coding to input the bibitems directly in the
% %% TeX file.

% %% Refer following link for more details about bibliography and citations.
% %% https://en.wikibooks.org/wiki/LaTeX/Bibliography_Management


% %% For authoryear reference style

% \bibitem[Lamport(1994)]{lamport94}
%   \textit{\LaTeX: a document preparation system},
%   Addison Wesley, Massachusetts,

% \end{thebibliography}
\end{document}